\documentclass{gtpart}     
\usepackage[english]{babel}
\usepackage{amssymb, amsfonts, amsmath, amsthm, amscd}
\usepackage{hyperref}
\usepackage{graphicx}
\usepackage{mathrsfs} 
\usepackage[all]{xy}
\usepackage{pinlabel}
\usepackage{color}

\title{Pair of pants decomposition of 4--manifolds} 

\author{Marco Golla}
\givenname{Marco}
\surname{Golla}
\address{Department of Mathematics, Box 480, 751 06 Uppsala, Sweden}
\email{marco dot golla at math dot uu dot se}

\author{Bruno Martelli}
\givenname{Bruno}
\surname{Martelli}
\address{Dipartimento di Matematica, Largo Pontecorvo 5, 56127 Pisa, Italy}
\email{martelli at dm dot unipi dot it}



\arxivreference{arXiv:1510.06325}
\arxivpassword{kzs2n}

\newtheorem{teo_intro}{Theorem}[]
\newtheorem{quest_intro}[teo_intro]{Question}

\newtheorem{lemma}{Lemma}[section]
\newtheorem{teo}[lemma]{Theorem}
\newtheorem{prop}[lemma]{Proposition}
\newtheorem{cor}[lemma]{Corollary}

\theoremstyle{definition}
\newtheorem{defn}[lemma]{Definition}

\newtheorem{probl}[lemma]{Problem}
\newtheorem{example}[lemma]{Example}

\theoremstyle{remark}
\newtheorem{oss}[lemma]{Remark} 
\newtheorem{rem}[lemma]{Remark}

\newcommand{\matR} {\ensuremath {\mathbb{R}}}

\newcommand{\matZ} {\ensuremath {\mathbb{Z}}}
\newcommand{\matC} {\ensuremath {\mathbb{C}}}
\newcommand{\matc}{\matC}

\newcommand{\matCP} {\ensuremath {\mathbb{CP}}}

\newcommand{\CPbar} {\ensuremath {\overline\matCP\vphantom\matCP^2}}

\newcommand{\matr} [4] {\left({\tiny\begin{array}{@{}c@{\ }c@{}} #1 & #2 \\ #3 & #4 \\ \end{array}} \right)}

\newcommand{\vect} [2] {\left({\tiny\begin{array}{@{}c@{\ }c@{}} #1 \\ #2 \\ \end{array}} \!\right)}

\newcommand{\nota} [1] {\caption{\footnotesize{#1}}}

\begin{document}

\begin{abstract}    
Using tropical geometry, Mikhalkin has proved that every smooth complex hypersurface in $\matCP^{n+1}$ decomposes into \emph{pairs of pants}: a pair of pants is a real compact $2n$--manifold with cornered boundary obtained by removing an open regular neighborhood of $n+2$ generic complex hyperplanes from $\matCP^n$.

As is well-known, every compact surface of genus $g\geqslant 2$ decomposes into pairs of pants, and it is now natural to investigate this construction in dimension 4. Which smooth closed 4--manifolds decompose into pair of pants? We address this problem here and construct many examples: we prove in particular that every finitely presented group is the fundamental group of a 4--manifold that decomposes into pairs of pants.
\end{abstract}

\maketitle

\section*{Introduction}
The decomposition of surfaces into pairs of pants is an extraordinary instrument in geometric topology that furnishes, among many other things, a nice parametrization for Teichm\"uller spaces. Mikhalkin has generalized this notion in \cite{M} to all even dimensions as follows: he defines the \emph{$2n$--dimensional pair of pants} as the manifold obtained by removing $n+2$ generic hyperplanes from $\matCP^n$. One actually removes open regular neighborhoods of the hyperplanes to get a compact real $2n$--manifold with stratified cornered boundary: when $n=1$ we get $\matCP^1$ minus three points, whence the usual pair of pants. 

Using some beautiful arguments from tropical geometry, Mikhalkin has proved in \cite{M} that every smooth complex hypersurface in $\matCP^{n+1}$ decomposes into pairs of pants. We address here the following natural question:

\begin{quest_intro} \label{main:quest}
Which smooth closed manifolds decompose into pairs of pants?
\end{quest_intro}

The question makes sense of course only for real smooth manifolds of even dimension $2n$. It is natural to expect the existence of many smooth manifolds that decompose into pairs of pants and are not complex projective hypersurfaces:
for instance, the hypersurfaces of $\matCP^2$ are precisely the closed orientable surfaces of genus $g=\frac {(d-1)(d-2)}2$ for some $d>0$, while by assembling pairs of pants we obtain closed orientable surfaces of any genus. 

In this paper we study pants decompositions in (real) dimension 4. We start by constructing explicit pants decompositions for some simple classes of closed 4--manifolds: $S^4$, torus bundles over surfaces, circle bundles over 3--dimensional graph manifolds, toric manifolds, the simply connected manifolds $\#_k(S^2\times S^2)$ and $\#_k\matCP^2\#_h\CPbar$. Then we prove the following theorem, which shows that the 4--manifolds that decompose into pairs of pants form a quite large class:

\begin{teo_intro}\label{main:teo:intro}
Every finitely presented group is the fundamental group of a closed 4--manifold that decomposes into pairs of pants.
\end{teo_intro}

We get in particular plenty of non-K\"ahler and hence non-projective 4-manifolds.
We expose a more detailed account of these results in the remaining part of this introduction.

\subsection*{Pair of pants decompositions}
Mikhalkin's definition of a pair of pants decomposition is slightly more flexible than the usual one adopted for surfaces: the boundary of a pair of pants (of any dimension $2n$) is naturally stratified into circle fibrations and an appropriate collapse of these circles is allowed. With this language, the sphere $S^2$ has a pants decomposition with a single pair of pants where each boundary component is collapsed to a point.

\begin{figure}
\centering
\includegraphics[scale=.8]{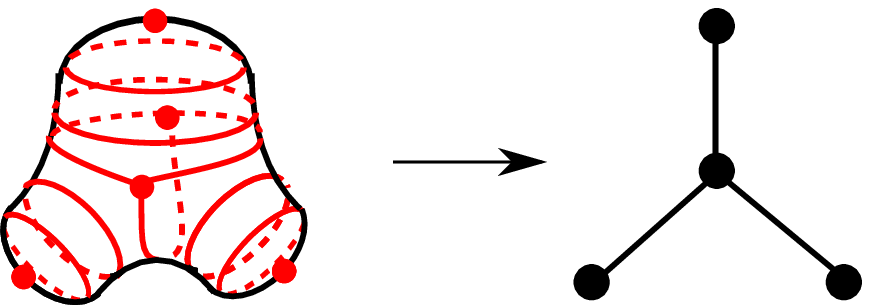}
\nota{The model fibration $\matCP^1 \to \Pi_1$.}\label{f:pantalone_filled:intro}
\end{figure}

\begin{figure}
\centering
\includegraphics[width = 12.5 cm]{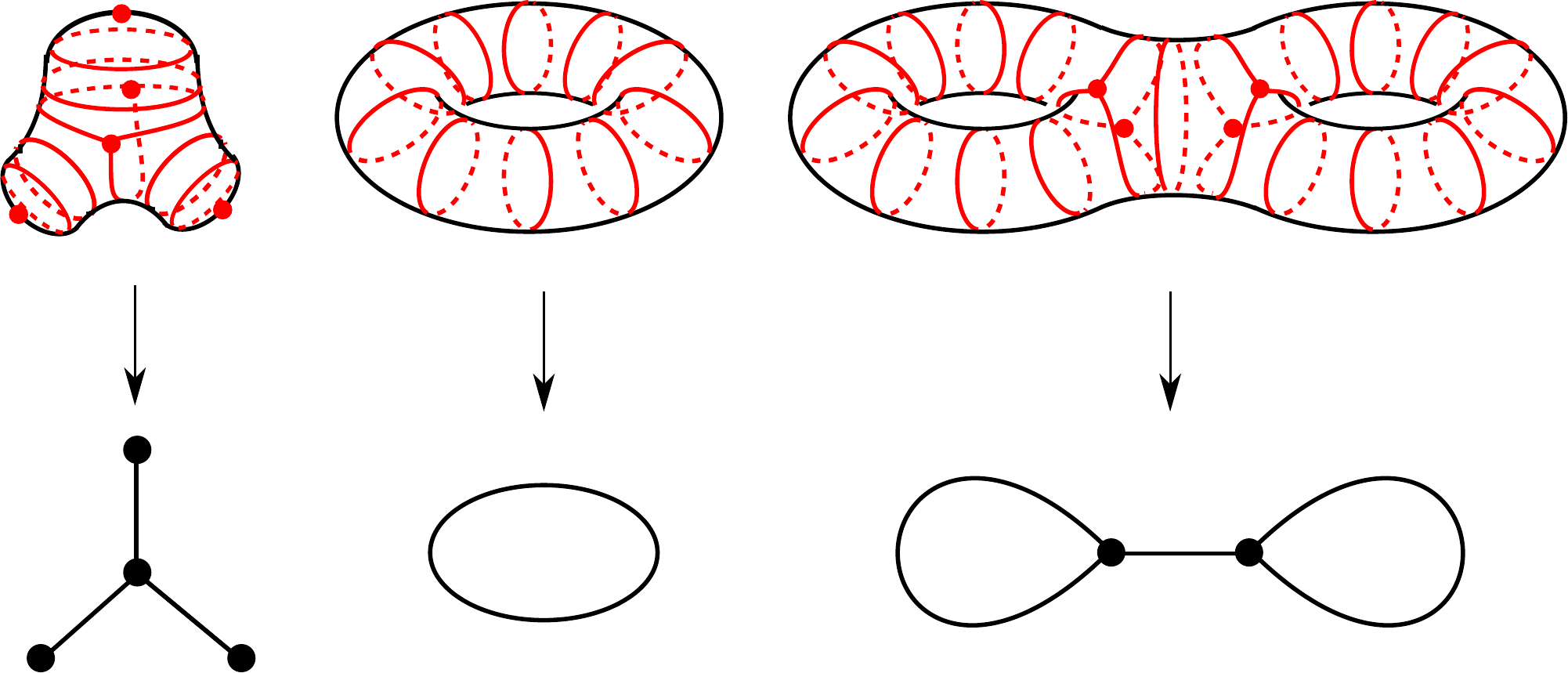}
\nota{Pair of pants decompositions of surfaces.}\label{f:examples}
\end{figure}

More precisely, a \emph{pair of pants decomposition} of a closed $2n$--manifold $M^{2n}$ is a fibration $M^{2n} \to X^n$ over a compact $n$--dimensional cell complex $X^n$ which is locally diffeomorphic to a \emph{model fibration} $\matCP^n \to \Pi_n$ derived from tropical geometry. The fiber of a generic (smooth) point in $X^n$ is a real $n$--torus.

When $n=1$, the model fibration is drawn in Figure
\ref{f:pantalone_filled:intro}
and some examples of pair of pants decompositions are shown in Figure \ref{f:examples}. The reader is invited to look at these pictures that, although quite elementary, describe some phenomena that will be present also in higher dimensions. When $n=1$, the base cell complex $X^1$ may be of these limited types: either a circle, or a graph with vertices of valence 1 and 3; there are three types of points $x$ in $X^1$ (smooth, a vertex with valence 1, or a vertex with valence 3), and the fiber over $x$ depends only on its type (a circle, a point, or a $\theta$-shaped graph, respectively).

In dimension $2n$, the model cell complex $\Pi_n$ is homeomorphic to the cone of the $(n-1)$--skeleton of the $(n+1)$--simplex, the model fibration sends $n+2$ generic hyperplanes onto the base $\partial \Pi_n$ of the cone and the complementary pair of pants onto its interior $\Pi_n \setminus \partial \Pi_n$. We are interested here in the case $n=2$.

\subsection*{Dimension 4}
In dimension 4, a pair of pants is $\matCP^2$ minus (the open regular neighborhood of) four generic lines: it is a 4--dimensional compact manifold with cornered boundary; the boundary consists of six copies of $P\times S^1$, where $P$ is the usual 2--dimensional pair-of-pants, bent along six 2--dimensional tori.

\begin{figure}
\labellist
\small\hair 2pt
\pinlabel $F_2$ at 41 100
\endlabellist
\centering
\includegraphics[scale=1.5]{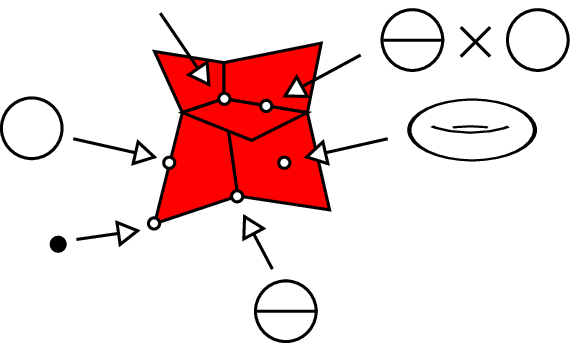}  
\nota{Fibers of the model fibration $\matCP^2 \to \Pi_2$.}\label{f:spine_fibers:intro}
\end{figure}

The model fibration $\matCP^2 \to \Pi_2$ is sketched in Figure \ref{f:spine_fibers:intro}: the cell complex $\Pi_2$ is homeomorphic to the cone over the 1--skeleton of a tetrahedron, and there are 5 types of points in $\Pi_2$; the fiber of a generic (i.e.~smooth) point of $\Pi_2$ is a torus, and the fibers over the other 4 types are: a point, $S^1$, $\theta$, $\theta \times S^1$, and a more complicated 2--dimensional
cell complex $F_2$ fibering over the central vertex of $\Pi_2$. The central fiber $F_2$ is a spine of the 4--dimensional pair of pants and is homotopically equivalent to a punctured 3-torus. (Likewise, when $n=1$ the fiber of the central vertex in $\matCP^1 \to \Pi_1$ is a $\theta$-shaped spine of the 2--dimensional 
pair of pants and is homotopically equivalent to a punctured 2--torus.)

A pair of pants decomposition of a closed 4--manifold $M^4$ is a map $M^4 \to X^2$ locally diffeomorphic to this model. The cell complex $X^2$ is locally diffeomorphic to $\Pi_2$, and the fiber over a generic point of $X^2$ is a torus.  


We note that pants decompositions are similar to (but different from) Turaev's \emph{shadows} \cite{Tu, CoThu}, that are fibrations $M^4 \to X^2$ onto similar cell complexes where the generic fiber is a disc, and $M^4$ is a 4--manifold with boundary that collapses onto $X^2$.

The main object of this work is to introduce many examples of 4--manifolds that decompose into pairs of pants. These examples are far from being exhaustive, and we are very far from having a satisfactory answer to Question \ref{main:quest}: for instance, we are not aware of any obstruction to the existence of a pants decomposition, so that we do not know if there is a closed 4--manifold which does not admit one.

\subsection*{Sketch of the proof of Theorem {\ref{main:teo:intro}}}
Once we set up the general theory, Theorem~\ref{main:teo:intro} is proved as follows. We first solve the problem of determining all the possible fibrations $M^4 \to X^2$ on a given $X^2$ by introducing an appropriate system of \emph{labelings} on $X^2$. We note that the same $X^2$ may admit fibrations $M^4\to X^2$ of different kinds, sometimes with pairwise non--diffeomorphic total spaces $M^4$, and each such fibration is detected by some labeling on $X^2$. This combinatorial encoding is interesting in its own because it furnishes a complete presentation of all pants decompositions in dimension 4. 

We then use these labelings to construct a large class of complexes $X^2$ for which there exist fibrations $M^4\to X^2$ that induce isomorphisms on fundamental groups. Finally, we show that every finitely presented $G$ has a $X^2$ in this class with $\pi_1(X^2) = G$.

\subsection*{Structure of the paper}
We introduce pair of pants decompositions in all dimensions in Section \ref{definitions:section}, following and expanding from Mikhalkin \cite{M} and focusing mainly on the 4--dimensional case. In Section \ref{four:section} we construct some examples.

In Section \ref{polygons:section} we study in detail the simple case when $X$ is a polygon. In this case $M\to X$ looks roughly like the moment map on a toric manifold, and every fibration $M\to X$ is encoded by some labeling on $X$. We then extend these labelings to more general complexes $X$ in Section \ref{general:section}.

In Section \ref{fundamental:groups:section} we prove Theorem \ref{main:teo:intro}.

\subsection*{Acknowledgements}
We thank Filippo Callegaro and Giovanni Gaiffi for stimulating discussions, and the referee for helpful comments and for finding a mistake in an earlier proof of Theorem~\ref{main:teo:intro}.

\section{Definitions} \label{definitions:section}

We introduce here simple complexes, tropical fibrations, and pair 
of pants decompositions. We describe these objects with some detail in dimension 2 and 4.

Recall that a \emph{subcomplex} $X\subset M$ of a smooth manifold $M$ is a subcomplex of some smooth triangulation of $M$. 

We work in the category of smooth manifolds: all the objects we consider are subcomplexes of some $\matR^N$, and a map between two such complexes is smooth if it locally extends to a smooth map on some open set.

\subsection{The basic cell complex $\Pi_n$}

\begin{figure}
\labellist
\small\hair 2pt
\pinlabel $[1,0,0]$ at 15 0
\pinlabel $[0,1,0]$ at 95 0
\pinlabel $[0,0,1]$ at 55 70
\pinlabel $[1,1,0]$ at 55 0
\pinlabel $[1,1,1]$ at 55 38
\pinlabel $[1,0,1]$ at 25 40
\pinlabel $[0,1,1]$ at 85 40
\endlabellist
\centering
\includegraphics[scale=1.5]{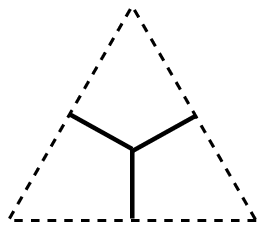}
\nota{The subcomplex $\Pi_1$ inside the standard simplex $\Delta$.}\label{f:Y}
\end{figure}

Let $\Delta$ be the standard $(n+1)$--simplex
\[\Delta = \big\{(x_0,\ldots, x_{n+1})\in \matR^{n+2} \ \big|\ x_0+\ldots + x_{n+1} = 1,\ x_i \geqslant 0\big\}.\]
We use the barycentric coordinates on $\Delta$, that is for every non-zero vector $x=(x_0, \ldots, x_{n+1}) \in \matR_{\geqslant 0}^{n+2}$
we denote by $[x_0, \ldots, x_{n+1}]$ the unique point in $\Delta$ that is a multiple of $x$. 
Every point $p \in \Delta$ has a unique description as $[x_0,\ldots, x_{n+1}]$ with $\max x_i = 1$ and we call it the \emph{normal form} of $p$.

\begin{figure}
\labellist
\small\hair 2pt
\pinlabel $(0,0)$ at 37 10
\pinlabel $(1,1)$ at 48 30
\pinlabel $(0,1)$ at 37 50
\pinlabel $(0,0)$ at 140 15
\pinlabel $(0,1)$ at 176 21
\pinlabel $(1,1)$ at 142 32
\pinlabel $(1,2)$ at 164 35
\pinlabel $(2,2)$ at 187 37
\pinlabel $(0,2)$ at 177 61
\endlabellist
\centering
\includegraphics[scale=1.5]{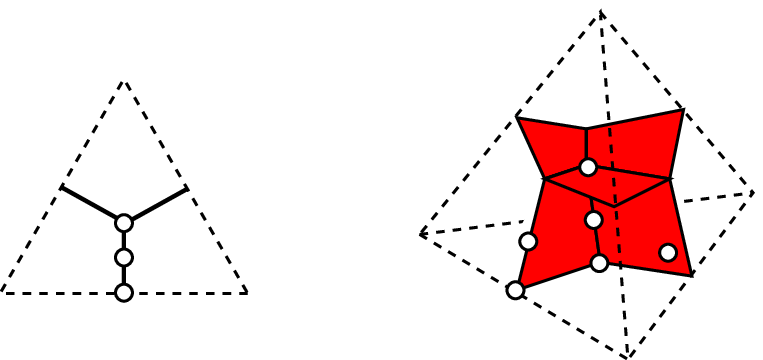}
\nota{The subcomplexes $\Pi_1$ and $\Pi_2$. Every point is of some type $(k,l)$ with $0\leqslant k\leqslant l\leqslant n$, and points of the same type define strata. Here $k$ is the dimension of the stratum and $l+1$ is the dimension of the face of $\Delta$ containing it.}\label{f:dualspine}
\end{figure}

\begin{defn}
Let $\Pi_n \subset \Delta$ be the following subcomplex:
\[\Pi_n = \big\{ [x_0, \ldots, x_{n+1}] \ |\ 0\leqslant x_i \leqslant 1 {\rm\ and\ } x_i = 1 {\rm\ for\ at\ least\ two\ values\ of\ } i \big\}.\]
\end{defn}
The subcomplex $\Pi_n \subset \Delta$ may be interpreted as the cut-locus of the vertices of $\Delta$, see $\Pi_1$ in Figure~\ref{f:Y}. Every point $x \in \Pi_n$ has a \emph{type} $(k,l)$ with $0\leqslant k \leqslant l \leqslant n$, which is determined by the following requirements: the normal form of $x$ contains $l-k+2$ different 1's and $n-l$ different 0's. More concretely, see Figure~\ref{f:dualspine} for the cases $n=1$ and $2$ which are of interest for us here. 

Points of the same type $(k,l)$ form some open $k$--cells, and these cells stratify $\Pi_n$. Geometrically, a point $x$ of type $(k,l)$ is contained in the $k$--stratum of $\Pi_n$ and in the $(l+1)$--stratum of $\Delta$. An open star neighborhood of $x$ in $\Pi_n$ is diffeomorphic to
the subcomplex
\[\Pi_{l,k} = \matR^k \times \Pi_{l-k} \times [0,+\infty)^{n-l}.\]
A point of type $(n,n)$ is a \emph{smooth} point, while the points with $l<n$ form the \emph{boundary} 
\[\partial \Pi_n = \Pi_n \cap \partial \Delta.\]

\subsection{Tropical fibration of $\matCP^n$}

\begin{figure}
\labellist
\small\hair 2pt
\pinlabel $[1,0,0]$ at 0 0
\pinlabel $[0,1,0]$ at 70 0
\pinlabel $[0,0,1]$ at 37 65
\pinlabel $[1,1,0]$ at 37 0
\pinlabel $[1,0,1]$ at 12 35
\pinlabel $[0,1,1]$ at 62 35
\endlabellist
\centering
\includegraphics[scale=1.8]{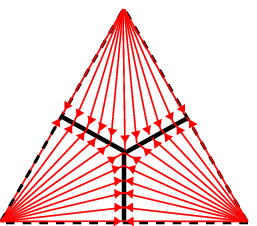}
\nota{The projection of $\Delta^*$ onto $\Pi_n$.}\label{f:Y_projection}
\end{figure}

Using tropical geometry, Mikhalkin has constructed in \cite{M} a map
\[\pi\colon\matCP^n \longrightarrow \Pi_n.\]
The map is defined as a composition $\pi = \pi_2 \circ \pi_1$ of two projections. The first one is a restriction of the projection
\begin{align*}
\matCP^{n+1} & \stackrel {\pi_1} \longrightarrow \Delta \\
[z_0, \ldots, z_{n+1}] & \longmapsto \big[|z_0|, \ldots, |z_{n+1}|\big].
\end{align*}
We identify $\matCP^n$ with the hyperplane $H\subset \matCP^{n+1}$ defined by the equation $z_0 + \ldots + z_{n+1}= 0$ and restrict $\pi_1$ to $H$. The image $\pi_1(H)$ is a region in $\Delta$ called \emph{amoeba} which contains $\Pi_n$ as a spine \cite{PR}. There is a simple projection that retracts the amoeba onto its spine $\Pi_n$: it is the restriction of a map
\[\pi_2\colon \Delta^* \longrightarrow \Pi_n\]
where $\Delta^*$ is $\Delta$ minus its vertices. The map $\pi_2$ is drawn in Figure~\ref{f:Y_projection} and is defined as follows: up to permuting the coordinates we suppose for simplicity that $x = [x_0, \ldots, x_{n+1}]$ with $x_0 \geqslant x_1 \geqslant \ldots \geqslant x_{n+1}$ and we define
\[\pi_2(x) = [x_1, x_1, x_2,\ldots, x_{n+1}].\]
The composition $\pi = \pi_2\circ \pi_1$ is a map that sends $\matCP^n = H$ onto $\Pi_n$.

The map $\pi_2$ is only piecewise smooth: it can then be smoothened as explained in \cite[Section 4.3]{M}, so that the composition $\pi$ is also smooth. In the following sections we study $\pi_2$ in the cases $n=1$ and $n=2$ before the smoothening, because it is easier to determine the fibers of $\pi$ concretely using the non-smoothed version of $\pi_2$. We remark that in the dimension 4 that we are interested in every piecewise-linear object can be easily smoothened, so this will not be an important issue anyway.

\subsection{The case $n=1$}
We now describe explicitly the fibration 
\[
\pi\colon \matCP^1 \to \Pi_1.
\] 
Recall that $\matCP^1$ is identified with the line $H= \{z_0 + z_1 + z_2 = 0\}$ in $\matCP^2$ and that $\Pi_1$ contains points of type $(0,0)$, $(1,1)$, and $(0,1)$. 

\begin{prop}
The fiber $\pi^{-1}(x)$ of a point $x\in\Pi_1$ is
\begin{itemize}
\item a point if $x$ is of type $(0,0)$,
\item a piecewise-smooth circle if $x$ is of type $(1,1)$,
\item a $\theta$-shaped smooth graph if $x$ is of type $(0,1)$.
\end{itemize}
\end{prop}

\begin{proof}
Up to reordering we have $x=[1,1,0]$, $[1,1,t]$, or $[1,1,1]$ with $0<t<1$ depending on the type. Using the calculation made in Figure~\ref{f:fiber_piecewise}-(left) we can describe the fibers explicitly:

\begin{align*}
\pi^{-1}\big([1,1,0]\big) & = [1,1,0], \\
\pi^{-1} \big([1,1,t] \big) & = \Big(\big\{[x,e^{i\theta},t]\ \big|\ |x|\geqslant 1\big\} 
\cup \big\{[e^{i\theta},x,t]\ \big|\ |x|\geqslant 1\big\}\Big) \cap H \\
& = \Big\{[-e^{i\theta} - t, e^{i\theta}, t]\ \Big|\ \cos\theta \geqslant -\frac t2 \Big\} \\
 & \quad  \cup \Big\{[e^{i\theta}, - e^{i\theta} - t, t]\ \Big|\ \cos\theta \geqslant -\frac t2 \Big\}, \\
\pi^{-1} \big([1,1,1] \big) & = \Big(\big\{[x,e^{i\theta},1]\ \big|\ |x|\geqslant 1\big\} 
\cup \big\{[e^{i\theta},1,x]\ \big|\ |x|\geqslant 1\big\} \\
& \quad \cup \big\{[1,x,e^{i\theta}]\ \big|\ |x|\geqslant 1\big\}
\Big) \cap H \\
& =
\Big\{[-e^{i\theta} - 1, e^{i\theta}, 1]\ \Big|\ \cos\theta \geqslant -\frac 12 \Big\} \\
 & \quad  \cup \Big\{[e^{i\theta}, 1, - e^{i\theta} - 1]\ \Big|\ \cos\theta \geqslant -\frac 12 \Big\} \\
  & \quad  \cup \Big\{[1, - e^{i\theta} - 1, e^{i\theta}]\ \Big|\ \cos\theta \geqslant -\frac 12 \Big\}.
\end{align*}

\begin{figure}
\labellist
\small\hair 2pt
\pinlabel $t$ at 120 105
\pinlabel $e^{i\theta}$ at 70 160
\pinlabel $-e^{i\theta}-t$ at 55 20
\pinlabel $u$ at 360 105
\pinlabel $e^{i\varphi}t$ at 300 130
\pinlabel $e^{i\varphi}t+u$ at 360 130
\pinlabel $e^{i\theta}$ at 260 130
\pinlabel $-e^{i\theta}-e^{i\varphi}t-u$ at 340 25
\endlabellist
\centering
\includegraphics[scale=.7]{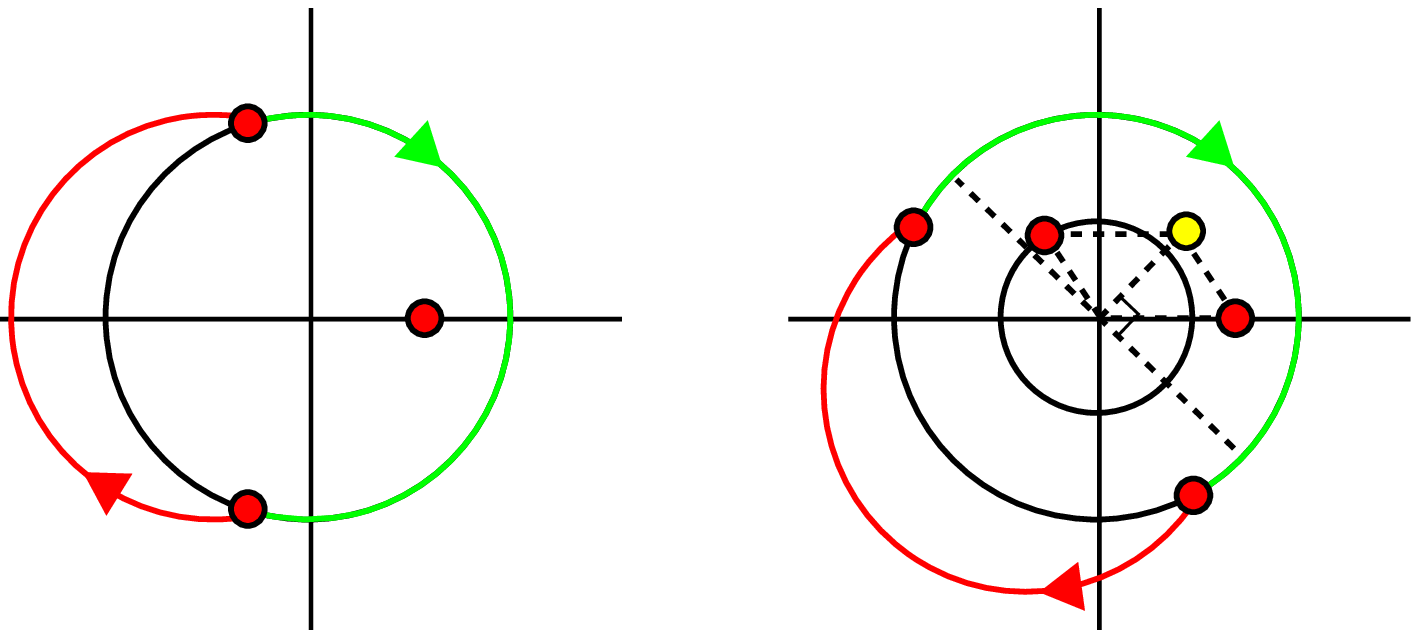}
\nota{The points $[-e^{i\theta}-t,e^{i\theta},t]$ with $|-e^{i\theta}-t|\geqslant 1$: when $\cos \theta = -\frac t2$ we get $-e^{i\theta} -t = e^{-i\theta} $; as the point $e^{i\theta}$ moves along the green arc of the unit circle the point $-e^{i\theta}-t$ moves along the red arc and has hence norm bigger than 1. This identifies one of the two arcs in $\pi^{-1}([1,1,t])$ (left). The fiber $\pi^{-1}([1,1,t,u])$ is considered similarly, with $e^{i\varphi}t + u$ instead of $t$ (right).}\label{f:fiber_piecewise}
\end{figure}

The fiber $\pi^{-1}([1,1,t])$ consists of two arcs with disjoint interiors but coinciding endpoints $[e^{\pm i\theta}, e^{\mp i\theta}, t]$ with $\cos\theta = -\frac t2$; therefore $\pi^{-1}([1,1,t])$ is a piecewise smooth circle. Analogously $\pi^{-1}([1,1,1])$ consists of three arcs joined at their endpoints $[e^{\pm \frac {2\pi i}3}, e^{\mp \frac {2\pi i}3}, 1]$ to form a $\theta$-shaped graph. 
\end{proof}

\begin{figure}
\centering
\includegraphics[scale=.8]{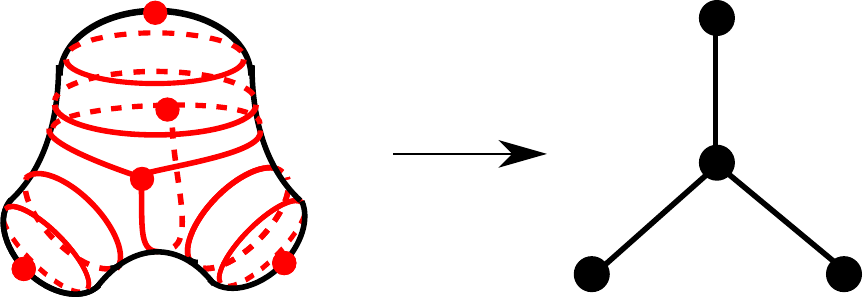}
\nota{The tropical fibration $\matCP^1 \to \Pi_1$.}\label{f:pantalone_filled}
\end{figure}

The fibration $\pi$ is homeomorphic to the one drawn in Figure~\ref{f:pantalone_filled}. The smoothing described in \cite[Section 4.3]{M} transforms the piecewise smooth circles into smooth circles, so that the resulting fibration is diffeomorphic to the one shown in the picture.

We note that the $\theta$-shaped graph is a spine of the pair of pants, and is also homotopic to a once-punctured 2--torus.
Both these facts generalize to higher dimensions.

\subsection{The case $n=2$}
We now study the fibration $\pi\colon \matCP^2 \to \Pi_2$, and our main goal is to show that its fibers are as in Figure \ref{f:spine_fibers}. 

Recall that we identify $\matCP^2$  with the plane $H = \{z_0+z_1+z_2+z_3 = 0\}$ in $\matCP^3$. The subcomplex $\Pi_2$ has points of type $(0,0), (1,1), (0,1)$ on its boundary, and of type $(2,2),(1,2),(0,2)$ in its interior. 

\begin{prop}
The fiber $\pi^{-1}(x)$ of a point $x\in \Pi_2$ is:
\begin{itemize}
\item a point if $x$ is of type $(0,0)$,
\item a piecewise-smooth circle if $x$ is of type $(1,1)$,
\item a $\theta$-shaped smooth graph $\theta$ if $x$ is of type $(0,1)$,
\item a piecewise-smooth torus if $x$ is of type $(2,2)$,
\item a piecewise-smooth product $\theta \times S^1$ if $x$ is of type $(1,2)$,
\item some 2--dimensional cell complex $F_2$ if $x$ is of type $(0,2)$.
\end{itemize}
\end{prop}
\begin{proof}
Up to reordering, the point $x$ is one of the following:
\[[1,1,0,0],\ [1,1,t,0],\ [1,1,1,0],\ [1,1,t,u],\ [1,1,1,t],\ [1,1,1,1]\]
with $1>t\geqslant u>0$.

\begin{figure}
\labellist
\small\hair 2pt
\pinlabel $F_2$ at 41 100
\endlabellist
\centering
\includegraphics[scale=1.5]{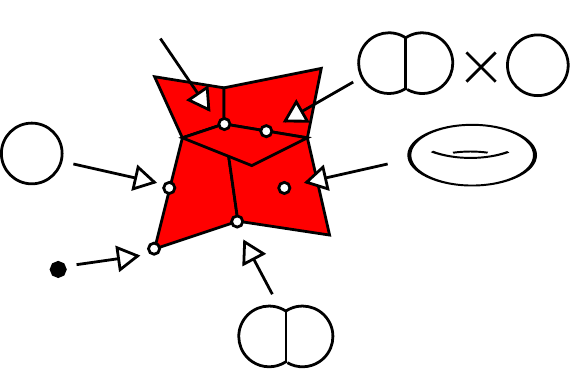}
\nota{Fibers of the tropical fibration $\matCP^2 \to \Pi_2$.}\label{f:spine_fibers}
\end{figure}

Let $f_1, f_2, f_3, f_4$ be the faces of $\Delta$, with $f_i = \{x_i = 0\}$. The preimage $\pi^{-1}_1 (f_i)$ is the plane $\{z_i = 0\}$ in $\matCP^3$ and intersects $H$ in a line $l_i$. The four lines $l_1,l_2,l_3,l_4$ are in general position in $H$ and intersect pairwise in the six points obtained by permuting the coordinates of $[1, -1, 0, 0]$.

The map $\pi$ sends $l_i$ onto the $Y$-shaped graph $f_i \cap \Pi_2$ exactly as described in the previous section, see Figure~\ref{f:pantalone_filled}. The map $\pi$ sends the four lines $l_i$ onto $\partial \Pi_2$, each line projected onto its own $Y$-shaped graph; the six intersection points are sent bijectively to the six points of type $(0,0)$ of $\Pi_2$.

It remains to understand the map $\pi$ over the interior of $\Pi_2$. Similarly as in the 1--dimensional case, Figure~\ref{f:fiber_piecewise}-(right) shows that

\begin{align*}
\pi^{-1} \big([1,1,t,u] \big) & = \Big(\big\{[x,e^{i\theta},e^{i\varphi}t,u]\big||x|\geqslant 1\big\} 
\cup \big\{[e^{i\theta},x,e^{i\varphi}t,u]\big||x|\geqslant 1\big\}\Big) \cap H \\
& = \Big\{[-e^{i\theta} - e^{i\varphi} t- u, e^{i\theta}, e^{i\varphi} t, u]\ \Big|\ \\
 & \qquad \cos (\theta - {\rm arg}(e^{i\varphi} t + u))  \geqslant -\frac 12 |e^{i\varphi} t + u|  \Big\} \\
 & \quad  \cup \big\{[e^{i\theta}, - e^{i\theta} - e^{i\varphi}t - u, e^{i\varphi}t, u]\ \Big|\ \\
 & \qquad \cos (\theta - {\rm arg}(e^{i\varphi} t + u))  \geqslant -\frac 12 |e^{i\varphi} t + u|  \Big\}.
 \end{align*} 

For every fixed $e^{i\varphi} \in S^1$ we get two arcs parametrized by $\theta$ with the same endpoints, thus forming a circle as in the 1--dimensional case. Therefore the fiber over $[1,1,t,u]$ is a (piecewise smooth) torus.

Analogously, the fiber over $[1,1,1,t]$ is a piecewise smooth product of a $\theta$-shaped graph and $S^1$. Finally, the fiber over $[1,1,1,1]$ is a more complicated 2--dimensional cell complex $F_2$. 
\end{proof}

The different fibers are shown in Figure~\ref{f:spine_fibers}. 
Let $F_i$ be the fiber over a point of type $(0,i)$. The fibers $F_0$, $F_1$, and $F_2$ are a point, a $\theta$-shaped graph, and some 2--dimensional complex.
These fibers ``generate'' all the others: the fiber over a point of type $(k,l)$ is piecewise-smoothly homeomorphic to $F_l \times (S^1)^k$. 


\begin{figure}
\labellist
\small\hair 2pt
\pinlabel $l_1$ at 170 19
\pinlabel $l_2$ at 95 145
\pinlabel $l_3$ at 75 45
\pinlabel $l_4$ at 35 95
\pinlabel $B$ at 180 140
\endlabellist
\centering
\includegraphics[scale=.8]{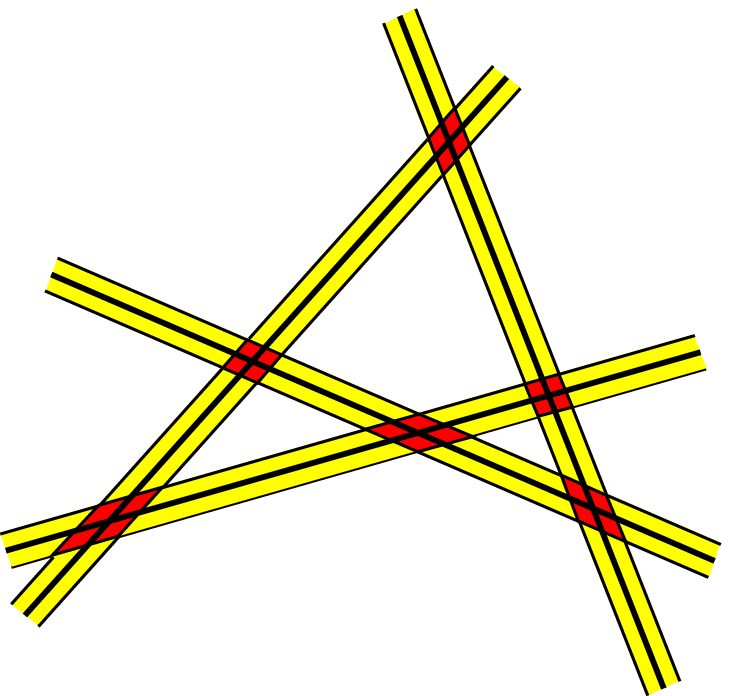}
\nota{A regular neighborhood of the four lines. It decomposes into six pieces diffeomorphic to $D^2 \times D^2$ (red) and four pieces diffeomorphic to $P\times D^2$ (yellow), where $P$ is a 2--dimensional pair of pants. (Every yellow piece is a $D^2$-bundle over $P$, and every such bundle is trivial. Note however that the normal bundle of each line is not trivial.) }\label{f:R}
\end{figure}

\begin{figure}
\centering
\includegraphics[scale=1.8]{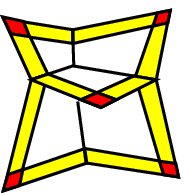}
\nota{A regular neighborhood of the four lines projects onto a regular neighborhood of $\partial \Pi_2$. Yellow and red blocks from Figure~\ref{f:R} project to the yellow and red portions in $\Pi_2$ drawn here. Note that there is a sixth sheet with a sixth red block behind the five that are shown.}\label{f:image_pi}
\end{figure}

\subsection{More on dimension 4} \label{more:4:subsection}
The fibration $\matCP^2 \to \Pi_2$ plays the main role in this work and we need to fully understand it. We consider here a couple of natural questions.

How does the fibration look like on a collar of $\partial \Pi_2$? It sends a regular neighborhood of the four lines $l_1, l_2, l_3, l_4$ shown in Figure~\ref{f:R} onto a regular neighborhood of $\partial \Pi_2$ drawn in Figure~\ref{f:image_pi}. Note that the regular neighborhood of the lines decomposes into pieces diffeomorphic to $D^2 \times D^2$ and $P\times D^2$, where $P$ is a pair of pants, see Figure~\ref{f:R}. On the red regions the fibration sends $D^2 \times D^2$ to $[0,1] \times [0,1]$ as $(w,z) \mapsto (|w|, |z|)$. On the yellow zone, each piece $P \times D^2$ is sent to $Y \times [0,1]$ as $(x,z) \mapsto (\pi(x), |z|)$ where $Y$ is a $Y$-shaped graph.

What is the fiber $F_2$? By construction it is a 2--dimensional spine of $\matCP^2$ minus the four lines. It is a well-known fact (proved for instance using the Salvetti complex \cite{S}) that the complement of four lines in general position in $\matCP^2$ is homotopically equivalent to a punctured 3--torus. More generally, the fiber $F_n$ is homotopic to a once-punctured $(n+1)$--torus (compare the case $n=1$). We have determined $F_2$ only up to homotopy, but this is sufficient for us.


\subsection{Simple complexes}
Always following Mikhalkin, we use the fibration $\pi$ as a standard model to define more general fibrations of manifolds onto complexes.

\begin{defn}
A \emph{simple $n$--dimensional complex} is a compact connected space $X \subset \matR^N$ such that every point has a neighborhood diffeomorphic to an open subset of $\Pi_n$. 
\end{defn}

For example, a simple 1--dimensional complex is either a circle or a graph with vertices of valence 1 and 3. 

Every point in $X$ inherits a type $(k,l)$ from $\Pi_n$, and points of the same type form a $k$--manifold called the \emph{$(k,l)$--stratum} of $X$. As opposed to $\Pi_n$, a connected component of a $(k,l)$--stratum needs not to be a cell: for instance, a closed smooth $n$--manifold is a simple complex where every point is smooth, i.e.~is of type $(n,n)$.

We use the word ``simple'' because it is largely employed to denote 2--dimensional complexes with generic singularities, see for instance \cite{Mat}. 

\subsection{Pants decomposition}
Let $M$ be a closed smooth manifold of dimension $2n$. Following \cite{M}, we define a \emph{pants decomposition} for $M$ to be a map
\[p\colon M \longrightarrow X\]
over a simple $n$--dimensional complex $X$ which is locally modeled on the fibration $\pi \colon \matCP^n \to \Pi_n$; that is, the following holds: for every point $x \in X$ there are an open neighborhood $U$ of $x$, a point $y$ in an open subset $V\subset\Pi_n$, a diffeomorphism $(U,x) \to (V,y)$, and a fiber-preserving diffeomorphism $\pi^{-1}(V) \to p^{-1}(U)$ such that the resulting diagram commutes:
\[ 
\xymatrix{ 
\pi^{-1} (V) \ar[r] \ar[d]_\pi & p^{-1}(U) \ar[d]^p  \\
V \ar[r] & U 
}
\]

\begin{figure}
\centering
\includegraphics[scale=.7]{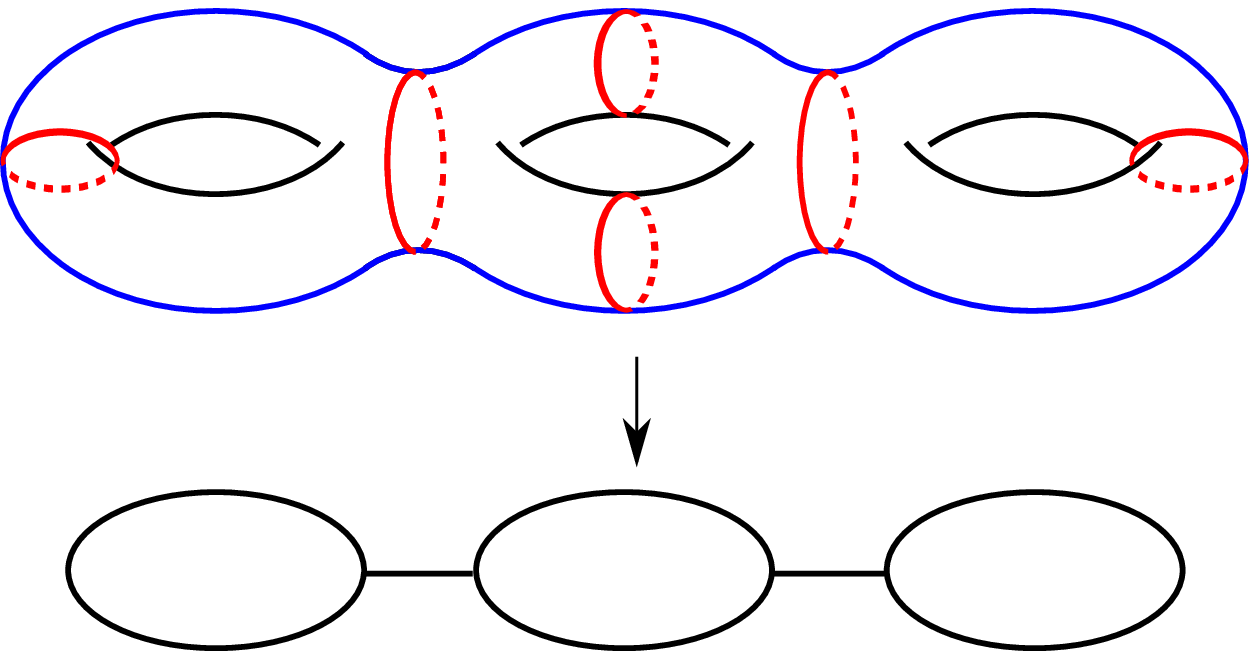}
\nota{A pants decomposition of a surface $M$ in the usual sense induces a fibration $M\to X$ onto a simple complex.}\label{f:pants_decomposition:fig}
\end{figure}

When $n=1$, a pants decomposition is a fibration $p\colon M\to X$ of a closed surface onto a 1--dimensional simple complex. If $X$ is not a circle and contains no 1--valent vertices, the fibration induces on $M$ a pants decomposition in the usual sense: the complex $X$ decomposes into $Y$-shaped subgraphs whose preimages in $M$ are pairs of pants, see Figure~\ref{f:pants_decomposition:fig}. Conversely, every usual pants decomposition of $M$ defines a fibration $M\to X$ of this type. 

In general the base complex $X$ may be quite flexible, for instance it might just be an $n$--manifold: therefore every smooth $n$--torus fibration on a $n$--manifold $X$ is a pants decomposition. Mikhalkin has proved the following remarkable result:

\begin{teo}[Mikhalkin \cite{M}]
Every smooth complex hypersurface in $\matCP^{n+1}$ admits a pants decomposition.
\end{teo}

As stated in the introduction, we would like to understand which manifolds of even dimension admit a pants decomposition. In dimension 2, 
every closed orientable surface has a pants decomposition: those having genus $g>1$ admit a usual one, while the sphere and the torus admit one in the more generalized sense introduced here; they fiber respectively over a segment (or a $Y$-shaped graph, or any tree) and a circle.

We now focus on the case $n=2$; that is, we look at smooth 4--manifolds fibering over simple 2--dimensional complexes.

\section{Four--manifolds} \label{four:section}

We now construct some closed 4--dimensional manifolds $M$ that decompose into pairs of pants, that is that admit a fibration $M\to X$ onto some simple complex $X$ locally modeled on $\matCP^2 \to \Pi_2$. In the subsequent sections we will study fibrations on a given $X$ in a more systematic way. 

\subsection{Some examples}
\begin{figure}
\centering
\includegraphics[scale=.7]{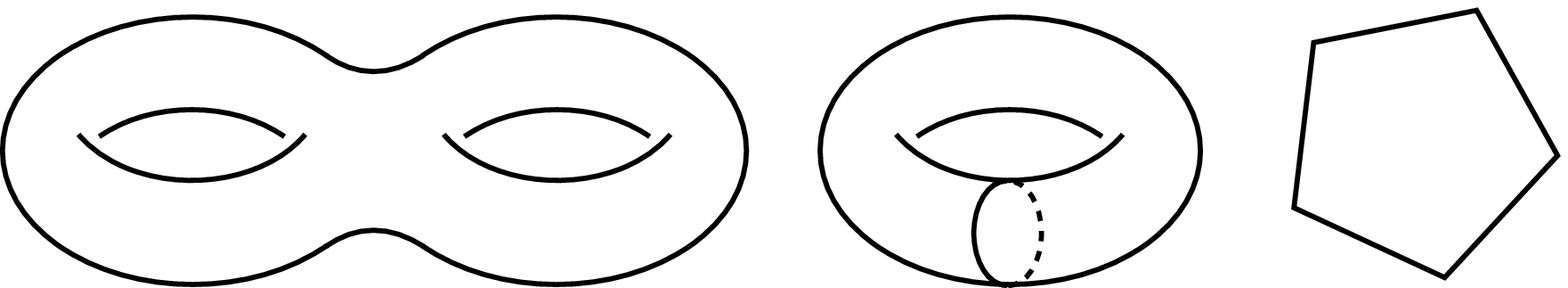}
\nota{Three simple 2--dimensional complexes: a closed surface (all points are smooth), a surface with a disc attached (all points are of type $(1,2)$ or $(2,2)$), and a polygon (all points are of type $(0,0)$, $(1,1)$, or $(2,2)$).}\label{f:simple_examples}
\end{figure}

We construct three families of examples of fibrations $M \to X$, that correspond to three simple types of complexes $X$ shown in Figure~\ref{f:simple_examples}: surfaces, surfaces with triple points, and polygons.

If $X$ is a closed surface, the fibrations $M \to X$ are precisely the torus bundles over $X$.

If $X$ contains points of type $(1,2)$ and $(2,2)$, we obtain more manifolds. Recall that a Waldhausen \emph{graph manifold} \cite{W} is any 3--manifold that decomposes along tori into pieces diffeomorphic to $P\times S^1$ and $D^2 \times S^1$, where $P$ is the pair of pants. For example, all lens spaces and Seifert manifolds are graph manifolds.

\begin{prop}
Let $p\colon M \to N$ be a circle bundle over an orientable closed graph manifold $N$. The closed manifold $M$ has a pants decomposition $M \to X$ for some $X$ consisting of points of type $(1,2)$ and $(2,2)$ only.
\end{prop}
\begin{proof}
It is proved in \cite[Proposition 3.31]{CoThu} that every orientable graph manifold $N$ admits a fibration $\pi$ 
over some simple complex $X$ called \emph{shadow} that consists of points of type $(1,2)$ and $(2,2)$ only, with fibers diffeomorphic respectively to a $\theta$-shaped graph and a circle. The composition of the two projections $\pi\circ p: M \to X$ is a pair of pants decomposition.
\end{proof}

If $X$ has only points of type $(0,0)$, $(1,1)$, and $(2,2)$, then it is a surface with polygonal boundary consisting of vertices and edges. Also in this case we get interesting manifolds.

\begin{prop}
A closed 4--dimensional toric manifold $M$ has a pants decomposition $M \to X$ for some polygonal disc $X$. In particular $\matCP^2$ fibers over the triangle.
\end{prop}
\begin{proof}
The moment map $M \to X$ is a fibration onto a polygon $X$ locally modeled on $\matCP^2 \to \Pi_2$ near a vertex of type $(0,0)$.
\end{proof}

The 4--dimensional 
closed toric varieties are $S^2 \times S^2$ and $\matCP^2 \# h \CPbar$, see \cite{FY}. In all the previous examples the base complex $X$ has no vertex of type $(0,2)$. 

\begin{probl}
Classify all the pair of pants decompositions $M \to X$ onto simple complexes $X$ without vertices of type $(0,2)$.
\end{probl}

This is a quite interesting set of not-too-complicated 4--manifolds, which contains torus bundles over surfaces, circle bundles over graph manifolds, and toric manifolds.

\subsection{Smooth hypersurfaces}
We now turn to more complicated examples where $X$ contains vertices of type $(0,2)$. Mikhalkin's theorem \cite[Theorem 1]{M} produces the following manifolds.

\begin{teo}
The smooth hypersurface $M$ of degree $d$ in $\matCP^3$ has a pants decomposition $M\to X$ on a simple complex $X$ with $d^3$ vertices of type $(0,2)$.
\end{teo}

Recall that the diffeomorphism type of $M$ depends only on the degree $d$. When $d=1,2,3,4$ the manifold $M$ is $\matCP^2$, $S^2 \times S^2$, $\matCP^2 \# 6\CPbar$, and the $K3$ surface, respectively.

\subsection{Euler characteristic}
The Euler characteristic of a pants decomposition can be easily calculated, and it depends only on the base $X$.
\begin{prop} \label{Euler:prop}
Let $M \to X$ be a pants decomposition. We have
\[\chi(M) = n_0 - n_1 + n_2\]
where $n_i$ is the number of points of type $(0,i)$ in $X$.
\end{prop}
\begin{proof}
All fibers have zero Euler characteristic, except the fibers $F_i$ above vertices of type $(0,i)$, that have $\chi(F_i) = (-1)^i$.
\end{proof}

\subsection{The nodal surface}
We note the following.
\begin{prop}
Let $p\colon M\to X$ be a pants decomposition. The preimage $S = p^{-1}(\partial X)$ is an immersed smooth compact surface in $M$.
\end{prop}
\begin{proof}
The fibration $p$ is locally modeled on the tropical fibration $\matCP^2 \to \Pi_2$ and the preimage of $\partial \Pi_2$ in $\matCP^2$ is an immersed surface consisting of four lines intersecting transversely in six points lying above the vertices of type $(0,0)$. 
\end{proof}
We call $S$ the \emph{nodal surface} of the fibration $p$.
It is an immersed surface in $M$ with one transverse self-intersection above each point of type $(0,0)$ of $X$. Every such self-intersection is called a \emph{node}.

\begin{rem}\label{r:nodalsphere}
We note that a vertex of type $(0,1)$ connected to three vertices of type $(0,0)$ determines an embedded sphere in $S$. Two vertices of type $(0,0)$ connected by an edge also determine an embedded sphere. 
\end{rem}

\section{Polygons} \label{polygons:section}

Let $X$ be a 2--dimensional simple complex. Is there a combinatorial way to encode all the pants decompositions $M\to X$ fibering over $X$? Yes, there is one, at least in the more restrictive case where every connected stratum in $X$ is a cell: every fibration is determined by some \emph{labeling} on $X$, which is roughly the assignment of some $2\times 2$ matrices to the connected 1--strata of $X$ satisfying some simple requirements. We describe this method here in the simple case when $X$ is a polygon. We will treat the general case in the next section.

\subsection{Fibrations over polygons}
\begin{figure}
\labellist
\small\hair 2pt
\pinlabel $S^1\times D^2$ at 270 15
\pinlabel $S^1\times S^1$ at 320 20
\pinlabel $D^2\times S^1$ at 322 80
\endlabellist
\centering
\includegraphics[scale=1]{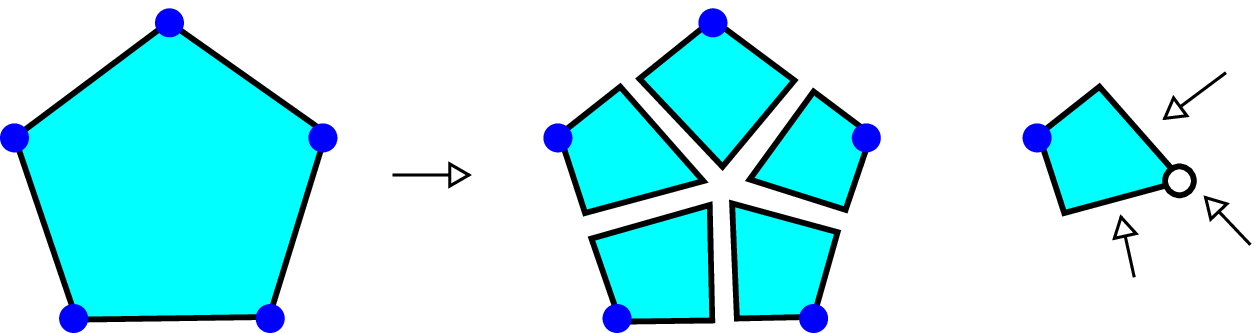}
\nota{A fibering $M\to X$ over a pentagon (left) can be broken into $n$ basic pieces (center). The fibering over each basi piece (right).}\label{f:polygon2}
\end{figure}

Let $X$ be a $n$--gon as in Figure \ref{f:polygon2}-(left), that is a simple 2--dimensional complex homeomorphic to a disc with $n\geqslant 1$ points of type $(0,0)$ called \emph{vertices}. The strata of type $(0,1)$ form $n$ \emph{edges} (or \emph{sides}).

Let $\pi\colon M \to X$ be a pair of pants decomposition. We first make some topological considerations.

\begin{prop}
The manifold $M$ is simply connected and $\chi(M) = n$. The nodal surface consists of $n$ spheres.
\end{prop}
\begin{proof}
We have $\chi(M)=n$ by Proposition \ref{Euler:prop}. The manifold $M$ is simply connected because $X$ is, and every loop contained in some fiber $\pi^{-1}(x)$ is homotopically trivial: it suffices to push $x$ to a vertex $v$ of $X$ and the loop contracts to the point $\pi^{-1}(v)$.

Thanks to Remark~\ref{r:nodalsphere}, the nodal surface consists of $n$ spheres, one above each edge of $X$.
\end{proof}

\subsection{Orientations}
In this paper we will be often concerned with orientations on manifolds, their products, and their boundaries. This can be an annoying source of mistakes, so we need to be careful. We will make use of the following formula on oriented manifolds $M$ and $N$:
\begin{equation} \label{ori:eqn}
\partial (M\times N) = (\partial M \times N) \cup (-1)^{\dim M} (M \times \partial N).
\end{equation}
Moreover, recall that the map
\begin{equation} \label{ori:i:eqn}
M\times N \longrightarrow N\times M
\end{equation}
that interchanges the two factors is orientation-preserving if and only if $\dim M \cdot \dim N$ is even.

\subsection{The basic fibration}
Let again $M\to X$ be a fibration over a polygon. We now break the given fibration $M\to X$ into some basic simple pieces, and show that $M\to X$ can be described by some simple combinatorial data.

\begin{figure}
\labellist
\small\hair 2pt
\pinlabel $L_1$ at 220 125
\pinlabel $L_2$ at 300 125
\pinlabel $L_3$ at 320 60
\pinlabel $L_4$ at 260 8
\pinlabel $L_5$ at 195 60
\endlabellist
\centering
\includegraphics[scale=0.8]{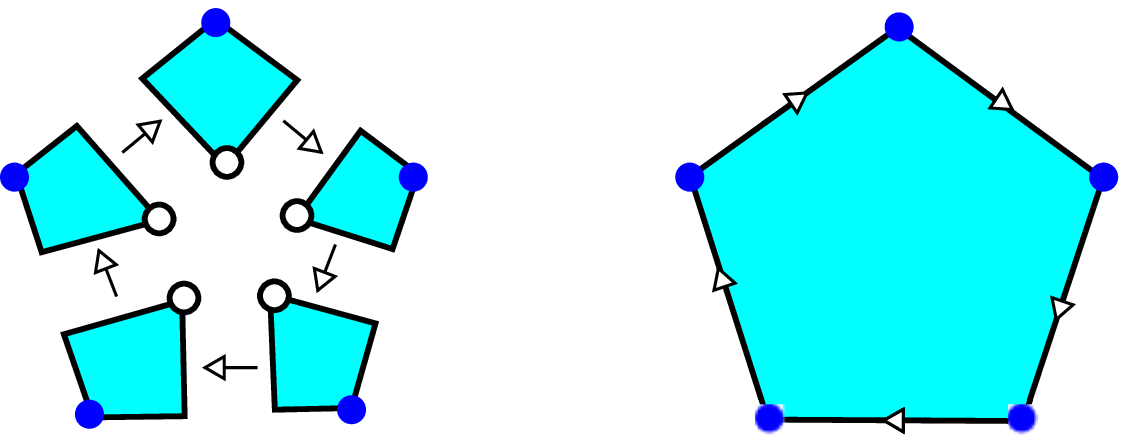}
\nota{The fibering $M\to X$ may be reconstructed by gluing the basic fibrations (left). The gluing can be determined by some label matrices $L_i$ (right).}\label{f:polygon}
\end{figure}

We break the $n$--gon into $n$ star neighborhoods of the vertices as in Figure \ref{f:polygon2}-(centre). Above each star neighborhood, the fibration is diffeomorphic to the \emph{basic fibration}
\[
D^2 \times D^2 \longrightarrow [0,1] \times [0,1]
\]
that sends $(w,z)$ to $(|w|,|z|)$, encountered in Section \ref{more:4:subsection} and sketched in Figure \ref{f:polygon2}-(right). The whole fibration $M\to X$ is constructed by gluing $n$ such basic fibrations as suggested in Figure \ref{f:polygon}-(left). We only need to find a combinatorial encoding of these gluings to determine $M\to X$.

Consider a single basic fibration $D^2\times D^2 \to [0,1]\times [0,1]$ as in Figure \ref{f:polygon2}-(right). The point $(0,0)$ is the fiber of $(0,0)$, the blue vertex in the figure. The boundary of $D^2 \times D^2$ is 
\[
\partial (D^2\times D^2) = (D^2\times S^1) \cup (S^1 \times D^2),
\] 
that is two solid tori (we call them \emph{facets}) cornered along the torus $S^1\times S^1$ (a \emph{ridge}). The manifold $D^2\times D^2$ is naturally oriented, and by (\ref{ori:eqn}) and (\ref{ori:i:eqn}) both solid tori inherit from $D^2\times D^2$ their natural orientations, which is invariant if we swap the factors $D^2$ and $S^1$. The ridge torus $S^1\times S^1$ however inherits opposite orientations from the two facets.

The ridge torus $S^1\times S^1$ is the fiber of $(1,1)$, the white dot in the figure, and the two facet solid tori fiber over the two adjacent sides $\{1\} \times [0,1]$ and $[0,1] \times \{1\}$.

Every arrow in Figure \ref{f:polygon}-(left) indicates a diffeomorphism $\psi\colon D^2 \times S^1 \to S^1 \times D^2 $ between two facets of two consecutive basic fibrations. It is convenient to write $\psi$ as a composition
\[
D^2\times S^1 \stackrel {\psi'} \longrightarrow D^2\times S^1 \stackrel j\longrightarrow S^1\times D^2
\]
where $j$ simply interchanges the two factors. By standard 3--manifolds theory, the diffeomorphism $\psi'$ is determined (up to isotopy) by its restriction to the boundary torus $S^1 \times S^1$, which is in turn determined (up to isotopy) by the integer invertible matrix $L\in {\rm GL}(2,\matZ)$ that encodes its action on $H_1(S^1\times S^1) = \matZ \times \matZ$. The only requirement is that $L$ must preserve the meridians, that is it must send $(1,0)$ to $(\pm 1,0)$. Summing up, we have the following.

\begin{prop}
The isotopy class of $\psi'$ is determined by a matrix 
\[
L = \begin{pmatrix} \varepsilon & k \\ 0 & \varepsilon' \end{pmatrix}
\]
with $\varepsilon, \varepsilon'  = \pm 1$ and $k\in \matZ$.
\end{prop}

We can encode all the gluings by assigning labels $L_1, \ldots, L_n$ of this type to the $n$ oriented edges of $X$ as in Figure \ref{f:polygon}-(right). We call such an assignment a \emph{labeling} of the polygon $X$. We define the matrices 
\[
I = \begin{pmatrix} 1 & 0 \\ 0 & 1 \end{pmatrix}, \qquad 
J=\begin{pmatrix} 0 & 1 \\ 1 & 0 \end{pmatrix}.
\]
Not every labeling defines a fibration $M\to X$. A necessary and sufficient condition is that the global monodromy around the central torus must be trivial.

\begin{prop}
The labeling defines a fibration $M\to X$ if and only if 
\[
JL_nJL_{n-1} \cdots JL_1 = I.
\]
If $\det L_i = -1$ for all $i$, the manifold $M$ is oriented.
\end{prop}
\begin{proof}
We only need to ensure that the monodromy around the central torus $S^1 \times S^1$ is isotopic to the identity, that is $JL_n \cdots JL_1 = I$. (The composition $\psi = \psi' \circ j$ translates into $JL$.) If $\det L_i=-1$ the standard orientations of the pieces $D^2\times D^2$ match to induce an orientation for $M$. \end{proof}

We say that the labeling is \emph{admissible} if $L_nJ\cdots L_1J = I$ and \emph{oriented} if $\det L_i = -1$ for all $i$. Summing up, we have proved the following.

\begin{prop}
Every fibration $M \to X$ over an $n$--gon $X$ is obtained by some admissible labeling on $X$.
\end{prop}

\begin{figure}
\labellist
\small\hair 2pt
\pinlabel $L$ at 70 115
\pinlabel $\matr 100{-1}$ at 218 130
\pinlabel $\matr {-1}001$ at 220 10
\pinlabel $\matr 1{-1}0{-1}$ at 320 90
\pinlabel $\matr 1{-1}0{-1}$ at 435 90
\pinlabel $\matr 1{-1}0{-1}$ at 380 0
\pinlabel $\matr 100{-1}$ at 490 140
\pinlabel $\matr 100{-1}$ at 610 140
\pinlabel $\matr 100{-1}$ at 610 20
\pinlabel $\matr 100{-1}$ at 490 20
\endlabellist
\centering
\includegraphics[scale=0.57]{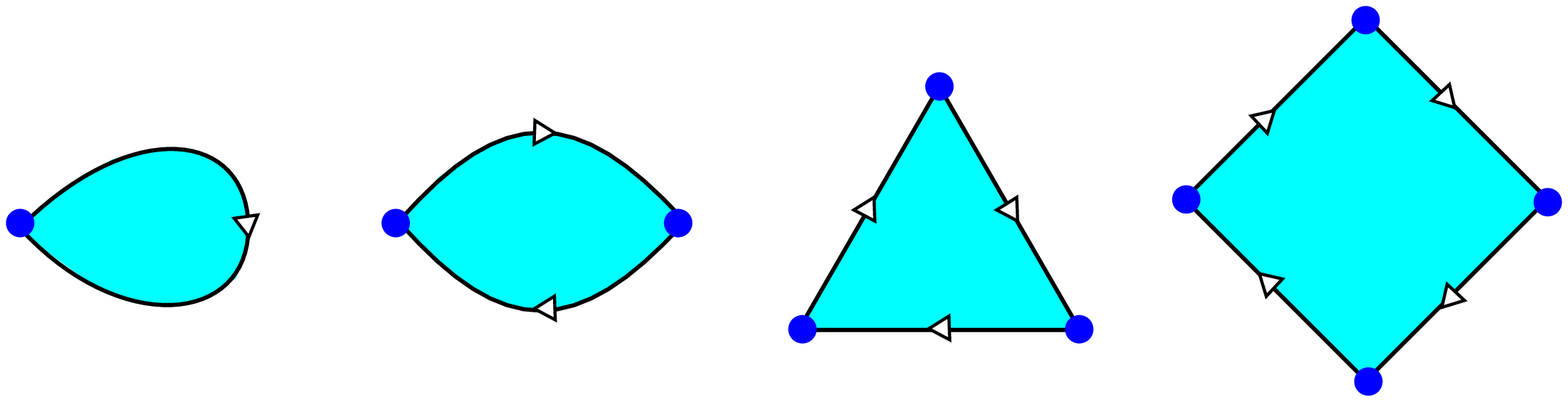}
\nota{The monogon has no admissible labeling. The other admissible labelings shown here represent $S^4$, $\matCP^2$, and $S^2\times S^2$.}\label{f:gons}
\end{figure}

Some examples are shown in Figure \ref{f:gons}. The monogon in Figure \ref{f:gons}-(left) has no admissible labeling $L$, because $LJ \neq I$ for every $L=\matr \varepsilon k 0 {\varepsilon'}$. The figure shows some oriented admissible labelings on the bigon, the triangle, and the square (admissibility is easily checked). Each determines a fibration $M\to X$.

\begin{prop}
The bigon in Figure \ref{f:gons} represents $S^4$.
\end{prop}
\begin{proof}
The bigon $X$ decomposes into two pieces $[0,1]\times [0,1]$ and $M$ decomposes correspondingly into two pieces $D^2 \times D^2$. The manifold $M$ decomposes smoothly into two 4--discs and is diffeomorphic to $S^4$.
\end{proof}

We have discovered that $S^4$ decomposes into pairs of pants. We will soon prove that the triangle and square in Figure \ref{f:gons} represent $\matCP^2$ and $S^2\times S^2$ respectively. 

Recall that we work entirely in the smooth (or equivalently, piecewise-linear) category. 

\subsection{Moves} \label{moves:subsection}
We now introduce some moves on admissible labelings.

Let $L_1,\ldots, L_n$ be a fixed oriented admissible labeling on the $n$--gon $X$ with edges $e_1,\ldots, e_n$. We know that it determines an oriented fibration $\pi\colon M\to X$. We start by noting that different labelings may yield the same fibration.

\begin{figure}
\labellist
\small\hair 2pt
\pinlabel $L_1$ at 15 30
\pinlabel $L_2$ at 12 170
\pinlabel $-L_1$ at 225 30
\pinlabel $-L_2$ at 220 170
\endlabellist
\centering
\includegraphics[scale=0.47]{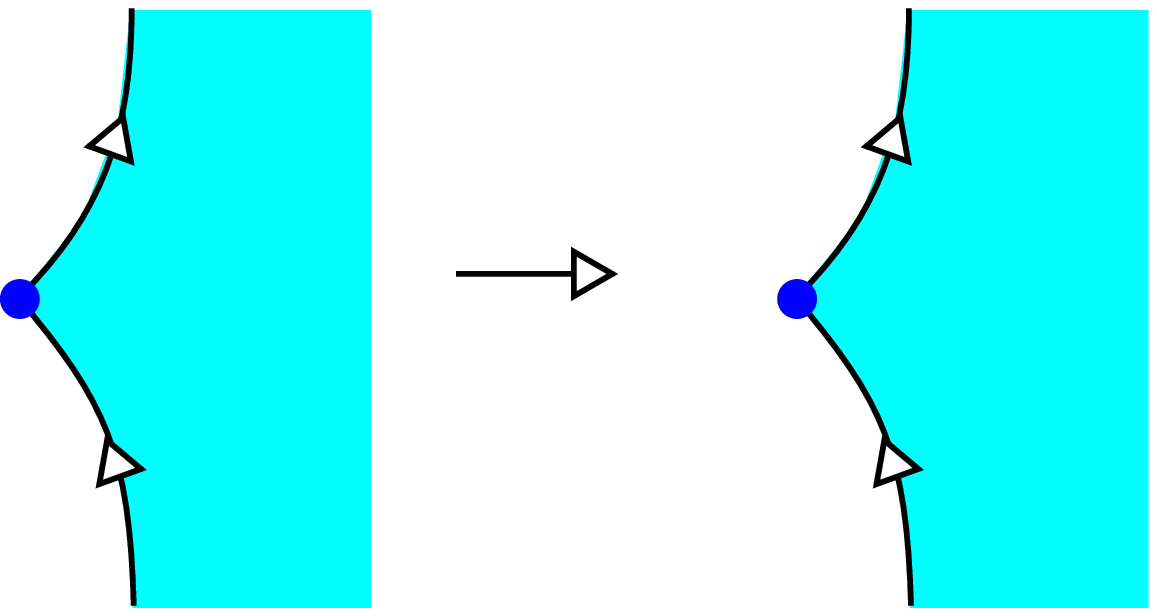}
\nota{If we change the signs of the labels of two consecutive edges, the fibration $M\to X$ remains unaffected.}\label{f:moves0}
\end{figure}

\begin{prop}\label{p:sign-move}
The move in Figure \ref{f:moves0} produces a new oriented admissible labeling, that encodes the same fibration $M\to X$.
\end{prop}
\begin{proof}
The fibration $D^2\times D^2 \to [0,1]\times [0,1]$ has the orientation-preserving automorphism $(z,w) \mapsto (\bar z, \bar w)$, that acts on $S^1\times S^1$ like $-I$. 
By employing it we see that the move produces isomorphic fibrations $M\to X$.
\end{proof}

Since $\det L_i = -1$ by hypothesis, every label $L_i$ is either $\matr 1{k_i}0{-1}$
or $\matr {-1}{k_i}0{1}$, and we say correspondingly that $L_i$ is \emph{positive} and \emph{negative}. By applying the move of Figure \ref{f:moves0} iteratively on the vertices of $X$ we may require that all labels $L_i$ are positive except at most one. Positive labels are preferable because of the following. 

Let $S_i$ be the sphere in the nodal surface lying above the edge $e_i$ of $X$. Note that two spheres $S_i$ and $S_j$ with $i\neq j$ intersect if and only if $e_i$ and $e_j$ are consecutive edges, and in that case they intersect transversely at the point (a node) projecting to the common vertex.

\begin{prop} \label{positive:prop}
If $L_i$ is positive, the sphere $S_i$ has a natural orientation. If $L_i$ and $L_{i+1}$ are positive, then $S_i\cdot S_{i+1} = +1$.
\end{prop}
\begin{proof}
The label $L_i$ represents the gluing of two pieces $D^2\times D^2$ and $D^2\times D^2$ along a map $\psi\colon D^2\times S^1 \to S^1\times D^2$ that sends the core $\{0\}\times S^1$ to $S^1\times \{0\}$. The sphere $S_i$ decomposes into two discs as $(\{0\} \times D^2) \cup_\psi (D^2\times \{0\})$.

If $L_i$ is positive, then $\psi$ identifies $\{0\}\times S^1$ to $ S^1\times \{0\}$ orientation-reversingly and hence the natural orientations of the two discs match to give an orientation for $S_i$.

The intersection of two consecutive $S_i$ and $S_j$ is transverse and positive (when they are both naturally oriented), because they intersect like $\{0\}\times D^2$ and $D^2\times \{0\}$ inside $D^2\times D^2$.
\end{proof}

Recall that the self-intersection $S_i\cdot S_i$ is independent of the chosen orientation for $S_i$ and is hence defined for all $i$, no matter whether $L_i$ is positive or not. The self-intersection of $S_i$ is easily detected by the labeling as follows.

\begin{prop} \label{pm:prop}
For each $i$, we have 
\[
L_i = \begin{pmatrix} {\pm 1} & {\mp (S_i\cdot S_i)} \\ 0 & {\mp 1} \end{pmatrix}
\]
for all $i$.
\end{prop}
\begin{proof}
Up to using the move in Figure \ref{f:moves0} we may restrict to the positive case $L_i = \matr 1k0{-1}$ and we need to prove that $S_i\cdot S_i = -k$. We calculate $S_i\cdot S_i$ by counting (with signs) the point in $S_i\cap S_i'$ where $S_i'$ is isotopic and transverse to $S_i$. 

Recall that $S_i = (\{0\}\times D^2) \cup_\psi (D^2 \times \{0\})$. We construct $S_i'$ by taking the discs $\{1\} \times D^2$ and $D^2 \times \{1\}$: their boundaries do not match in $S^1 \times S^1$ because they form two distinct longitudes in the boundary of the solid torus $S^1 \times D^2$, of type $(1,0)$ and $(1,k)$. We can isotope the former longitude to the latter inside the solid torus, at the price of intersecting the core $S^1 \times 0$ some $|k|$ times: in this way we get a $S_i'$ that intersects $S_i$ transversely into these $|k|$ points, always with the same sign.

We have proved that $S_i\cdot S_i = \pm k$. To determine the sign, it suffices to consider one specific case. We pick the triangle $X$ in Figure \ref{f:gons}, where all labels are $\matr 1{-1}0{-1}$. Here $\chi(M) = 3$ and $M$ is simply connected, therefore $H_2(M) = \matZ$. The nodal surface contains three spheres $S_1,S_2,S_3$ that represent  elements in $H_2(M)$ with $S_i\cdot S_i = \varepsilon = \pm 1$ for each $i$
 and $S_i \cdot S_j = 1$ when $i\neq j$. In particular, $S_i$ is a generator of $H_2(M)$ for each $i$.
Since $S_1\cdot S_2 = S_1\cdot S_3 = 1$, then $S_2 = S_3 = \varepsilon S_1$, and hence $1 = S_2\cdot S_3 = \varepsilon^2 \cdot S_1\cdot S_1 = \varepsilon$, hence $S_i\cdot S_i = +1$.
\end{proof}

\begin{figure}
\labellist
\small\hair 2pt
\pinlabel $L_1$ at 15 250
\pinlabel $L_2$ at 12 390
\pinlabel $L_1\cdot\matr {-1}10{-1}$ at 205 250
\pinlabel $\matr 1{-1}0{-1}$ at 210 320
\pinlabel $\matr 1101\cdot L_2$ at 210 390
\pinlabel $L_1$ at 412 250
\pinlabel $L_2$ at 409 390
\pinlabel $L_1\cdot\matr 1101$ at 602 250
\pinlabel $\matr 110{-1}$ at 607 320
\pinlabel $\matr 1{-1}01\cdot L_2$ at 603 390
\pinlabel $L$ at 445 30
\pinlabel $L$ at 205 100
\pinlabel $\matr {-1}001$ at 425 195
\pinlabel $\matr 100{-1}$ at 402 100
\endlabellist
\centering
\includegraphics[scale=0.47]{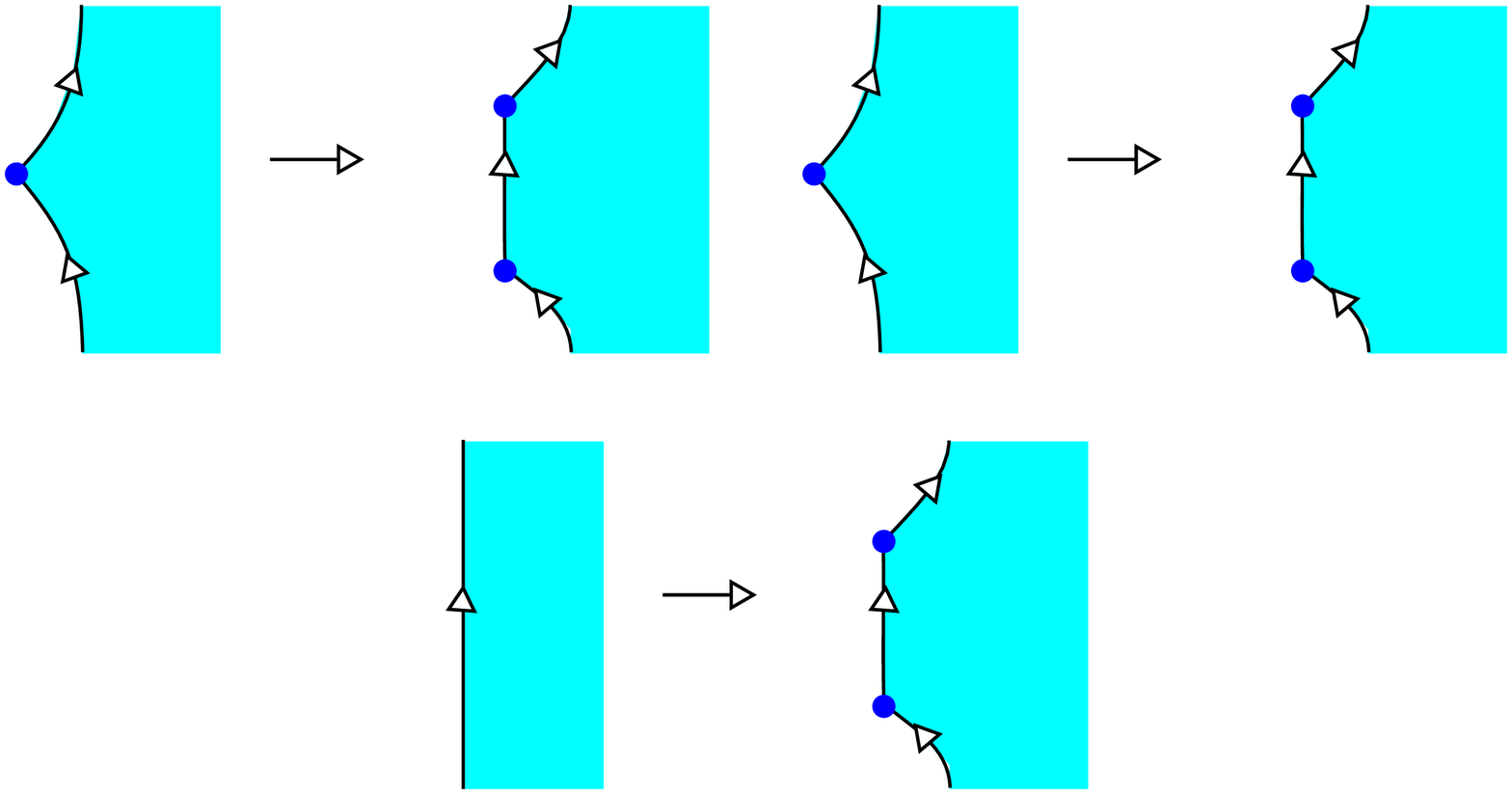}
\nota{Three moves that transform $M$ by connected sum with $\matCP^2$ (top left), $\overline \matCP^2$ (top right), and $S^2\times S^2$ (bottom centre).}\label{f:moves}
\end{figure}

We now consider two more moves on positive admissible labelings, shown in Figure \ref{f:moves}. It is easily checked that they both transform $L_1,\ldots, L_n$ into a new positive admissible labeling on a bigger polygon.

\begin{prop} \label{blow:prop}
The three moves in Figure \ref{f:moves} transform $M$ into 
\[
M\# \matCP^2, \qquad M\#\CPbar, \qquad M\#(S^2\times S^2)
\]
respectively.
\end{prop}

\begin{proof}
Both moves transform a fibration $M\to X$ into a new fibration $M'\to X'$.
The first two moves substitute a vertex $v$ of $X$ with a new edge $e$. The preimages of $v$ and $e$ in $M$ and $M'$ are a point $x\in M$ and a sphere $S\subset M'$ with $S\cdot S = +1$ or $S\cdot S=-1$ depending on the move. Substituting $x$ with $S$ amounts to making a topological blowup, that is a connected sum with $\matCP^2$ and $\CPbar$, respectively. 

The third move substitutes a point $x$ contained in some edge of $X$ with a new edge $e$. The preimages of $x$ and $e$ in $M$ and $M'$ are a circle $\gamma$ and a sphere $S\subset M'$ with $S\cdot S = 0$. The substitution of $\gamma$ with $S$ is called a \emph{surgery}, and since $M$ is simply connected the effect is a connected sum with $S^2\times S^2$.
\end{proof}
 
In particular the triangle and square from Figure \ref{f:gons} represent the oriented smooth 4--manifolds $\matCP^2$ and $S^2\times S^2$.

\begin{cor} \label{ngon:cor}
If $M= \#_h\matCP^2\#_k\CPbar$ or $M=\#_h (S^2\times S^2)$, then $M$ decomposes into pairs of pants; more precisely, $M$ fibers over the $n$--gon, with $n=\chi(M)$.
\end{cor}

These oriented manifolds are in fact all we can get from a polygon $X$.

\begin{prop}
Every oriented labeling on a polygon $X$ represents one of the manifolds of Corollary \ref{ngon:cor}.
\end{prop}
\begin{proof}
Every label is of type $L_i = \matr {\pm 1} {h_i} 0 {\mp 1}$. If $|h_i|\leqslant 1$, we can simplify $X$ via one of the moves from Figure \ref{f:moves} and proceed by induction. If $|h_i|\geqslant 2$ for all $i$, it is easy to show that the coloring cannot be admissible, because the product $L_nJ\cdots L_1J$ cannot be equal to $I$.

Indeed, we have $M_i = L_iJ = \matr{h_i}{\pm 1}{\mp 1}0$. The matrix $M_1$ sends $\vect 10$ to some $\vect ab$ with $|a|>|b|>0$, and any such vector is sent by any $M_i$ to a vector $\vect {a'}{b'}$ with $|a'|>|b'|>0$ again, so $M_n \cdots M_1 \vect 10 \neq \vect 10$.
\end{proof} 
 
\section{The general case} \label{general:section}
We now extend the discussion of the previous section from polygons to more general simple complexes $X$. For the sake of simplicity, we restrict our investigation to a class of complexes called \emph{special}, whose strata are all discs.
 
\begin{defn}\label{d:special}
A simple complex $X$ is \emph{special} if the connected components of all the $(k,l)$--strata are open $k$--cells.
\end{defn}

For instance, the polygons and
the model complex $\Pi_n$ are special. Every connected component of each stratum in a special 2--dimensional complex $X$ is a cell, called \emph{vertex}, \emph{edge}, or \emph{face} according to its dimension. Vertices are of type $(0,0)$, $(0,1)$, or $(0,2)$, and edges are of type $(1,1)$ or $(1,2)$. Each face is a polygon with $m$ edges and $m$ vertices for some $m$, and the vertices may be of different types.

\subsection{The basic fibrations}
Let $M\to X$ be a fibration over some special complex $X$. We now extend the discussion of the previous section to this more general setting: we break $M\to X$ into basic fibrations of three types, and we show that $M\to X$ may be encoded by some combinatorial labeling on $X$ that indicate the way these basic fibrations match along their (cornered) boundaries.

\begin{figure}
\centering
\includegraphics[scale=1.5]{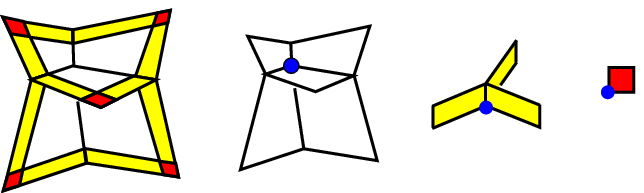}
\nota{The complex $\Pi_1$ (left) decomposes into the star neighborhoods of its vertices (right).}\label{f:blocks2}
\end{figure}

A $n$--gon breaks into $n$ star neighborhoods of its vertices as in Figure \ref{f:polygon2}; analogously, every special complex $X$ decomposes into star neighborhoods $S_v$ of its vertices $v$, which are now of three different types $(0,0)$, $(0,1)$, and $(0,2)$. 
For instance, the model complex $\Pi_2$ decomposes into 11 pieces, as shown in Figure~\ref{f:blocks2}: these are 6, 4, 1 stars of vertices of type $(0,0)$, $(0,1)$, $(0,2)$ respectively.

\begin{figure}
\labellist
\small\hair 2pt
\pinlabel $P\times S^1$ at 35 40
\pinlabel $P\times S^1$ at 71 84
\pinlabel $P\times S^1$ at 71 12
\pinlabel $P\times S^1$ at 107 66
\pinlabel $P\times S^1$ at 128 70
\pinlabel $D^2\times S^1$ at 105 42
\pinlabel $D^2\times S^1$ at 177 42
\pinlabel $D^2\times S^1$ at 163 70
\pinlabel $D^2\times S^1$ at 188 75
\pinlabel $D^2\times S^1$ at 215 52
\endlabellist
\centering
\includegraphics[scale=1.8]{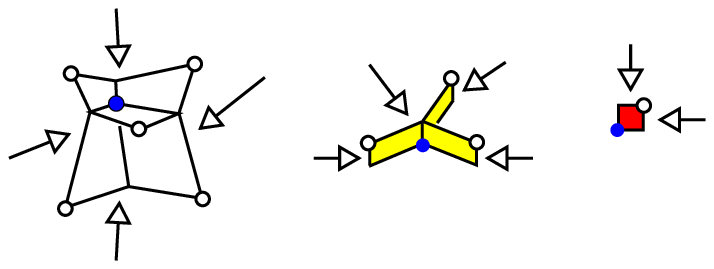}
\nota{Every block $M_v$ is a compact 4--manifold with corners: its boundary is a closed connected 3--manifold cornered along tori. There is one corner torus above each white dot, and the tori decompose the 3--manifold $\partial M_v$ into pieces diffeomorphic to $S^1\times P$ or $S^1 \times D^2$. Here $P$ indicates the pair of pants}\label{f:blocks3}
\end{figure}

The fibration $\pi\colon \matCP^2 \to \Pi_2$ decomposes correspondingly into $6+4+1=11$ 
\emph{basic fibrations} $M_v \to S_v$ above the star neighborhood $S_v$ of each vertex $v$. Every manifold $M_v$ is a regular neighborhood in $\matCP^2$ of the fiber $\pi^{-1}(v)$ of $v$, and its topology is deduced from Figures~\ref{f:R} and \ref{f:image_pi}.

There are three basic fibrations $M_v\to S_v$ to analyze, depending on the type of the vertex $v$. If $v$ is of type $(0,0)$ or $(0,1)$ the fibration $M_v \to S_v$ is diffeomorphic to the following:
\[
\begin{array}{rcl} 
D^2 \times D^2 & \longrightarrow & [0,1] \times [0,1] \\
(z,w) & \longmapsto & (|z|,|w|)
\end{array}
\qquad 
\begin{array}{rcl} 
D^2 \times P & \longrightarrow & [0,1] \times Y\\
(z,x) & \longmapsto & (|z|, \pi(x))
\end{array}
\]
where $Y$ is a $Y$-shaped graph and $\pi\colon P \to Y$ is the tropical fibration, see Figure \ref{f:blocks3}. 
Both $D^2$ and $P$ are naturally oriented as subsets of some complex line in $\matCP^2$.
In both cases $M_v$ is a product and its boundary is 
\begin{align*}
\partial(D^2 \times D^2) & =  (D^2\times S^1) \cup (S^1 \times D^2), \\
\partial(D^2 \times P) & = (D^2\times \partial P) \cup (S^1 \times P).
\end{align*}
Recall the orientation formulas (\ref{ori:eqn}) and (\ref{ori:i:eqn}). The boundary consists of some \emph{facets} cornered along tori (the \emph{ridges}). The facets are either solid tori or $S^1 \times P$. We identify once for all orientation-preservingly every boundary component of $P$ with $S^1$, so that $D^2\times \partial P$ is identified to three copies of $D^2\times S^1$. There are three corner tori in $S^1\times \partial P$. 

\subsection{The pair of pants}
If $v$ is of type $(0,2)$, the block $M_v$ is not a product: it is the compact \emph{pair of pants}, as named by Mikhalkin \cite{M}, diffeomorphic to the complement of an open regular neighborhood of four generic lines $l_1,\ldots, l_4$ in $\matCP^2$. Its boundary has four facets $f_1,\ldots, f_4$, each diffeomorphic to $S^1 \times P$, cornered along six tori, one for each pair $l_i, l_j$ of distinct lines.

The facet $f_i$ is an $S^1$--bundle over some pair of pants $P_i\subset l_i$ obtained from $l_i$ by removing open discs containing the intersection points with the other lines. The bundle is necessarily trivial, since it is a circle bundle over a compact orientable surface with non-empty boundary; hence $f_i$ is diffeomorphic to $S^1 \times P$, but unfortunately \emph{not} in a canonical way (not even up to isotopy): the diffeomorphism depends (up to isotopy) on the choice of a \emph{section} of the bundle, and on an orientation of the fibers (this is a standard fact on circle bundles over surfaces with boundary).

A natural way to construct a section goes as follows.
Pick a line $r \in \matCP^2$ that intersects $l_i$ in one of the three points $l_i\cap \l_j$, for some $j\neq i$. The line $r$ provides a section
of the normal bundle of $l_i$ that vanishes only at $l_i\cap l_j$, and hence a section of the circle bundle over $P_i$. The isotopy class of the section in fact does not depend on the chosen line $r$, but only on the point $l_i\cap l_j$, so there are three possible choices.

We now fix an arbitrary partition $\{l_1,l_2\}, \{l_3,l_4\}$ of the four lines into two pairs, and define $r$ to be the line passing through the points $l_1\cap l_2$ and $l_3\cap l_4$. We use the line $r$ to define sections on all the four facets $f_i$ simultaneously as just explained. 

Each section is oriented as a subset of $r$ and identified with $P$. To complete the identification of $f_i$ with $S^1 \times P$ we need to orient the fibers: we orient them so that $S^1 \times P$ gets the correct orientation as a boundary portion of the block $M_v$ (which is in turn oriented as a domain in $\matCP^2$).

\begin{oss}
By taking an affine chart that sends $r$ to infinity, we see that 
\[
\matCP^2\setminus (l_1\cup l_2 \cup l_3 \cup l_4 \cup r) \cong \big(\matC\setminus \{0,1\} \big) \times \big( \matC \setminus \{0,1\} \big)
\]
so $M_v$ minus an open neighborhood of $r$ is naturally diffeomorphic to a product $P\times P$. This diffeomorphism furnishes the identifications of each $f_i$ with $S^1 \times P$ just described.
\end{oss}

\begin{figure}
\labellist
\small\hair 2pt
\pinlabel ${\matr 1{-1}0{-1}}$ at 72 25
\pinlabel ${\matr 1{-1}0{-1}}$ at 96 3
\pinlabel ${\matr 1{-1}0{-1}}$ at 113 67
\pinlabel ${\matr 1{-1}0{-1}}$ at 138 50
\endlabellist
\centering
\includegraphics[scale=1.4]{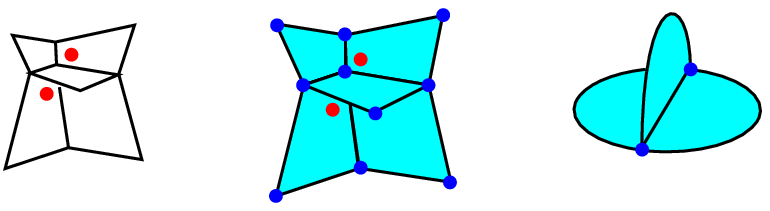}
\nota{At every vertex of type $(0,2)$, we fix two (of the six) opposite faces and we mark them with red dots (left). An admissible oriented labeling on $X$ that represents the tropical fibration $\matCP^2 \to \Pi_2$ (centre) and one that represents $(S^1\times S^3)\# (S^1\times S^3)$ (right).} \label{f:dots}
\end{figure}

There are of course three possible partitions of $\{l_1,l_2,l_3,l_4\}$ to choose from. To indicate on $X$ which partition we use, we mark with a dot the two opposite faces near $v$ that correspond to the pairs $l_1,l_2$ and $l_3,l_4$, as in Figure \ref{f:dots}-(left). This mark fixes an identification of every facet $f_i$ with the product $S^1 \times P$. 

\subsection{Labeling}
Every fibration $M\to X$ decomposes into basic fibrations, glued along facets that are either $D\times S^1$ or $S^1 \times P$. We now encode every such gluing with an appropriate labeling on $X$, that extends the one introduced in Section \ref{polygons:section} for polygons.

\begin{figure}
\labellist
\small\hair 2pt
\pinlabel $f$ at 200 130
\pinlabel $L_1$ at 100 30
\pinlabel $L_2$ at 80 150
\pinlabel $L_3$ at 180 230
\pinlabel $L_4$ at 290 230
\pinlabel $L_5$ at 340 150
\pinlabel $L_6$ at 240 30
\endlabellist
\centering
\includegraphics[scale=.45]{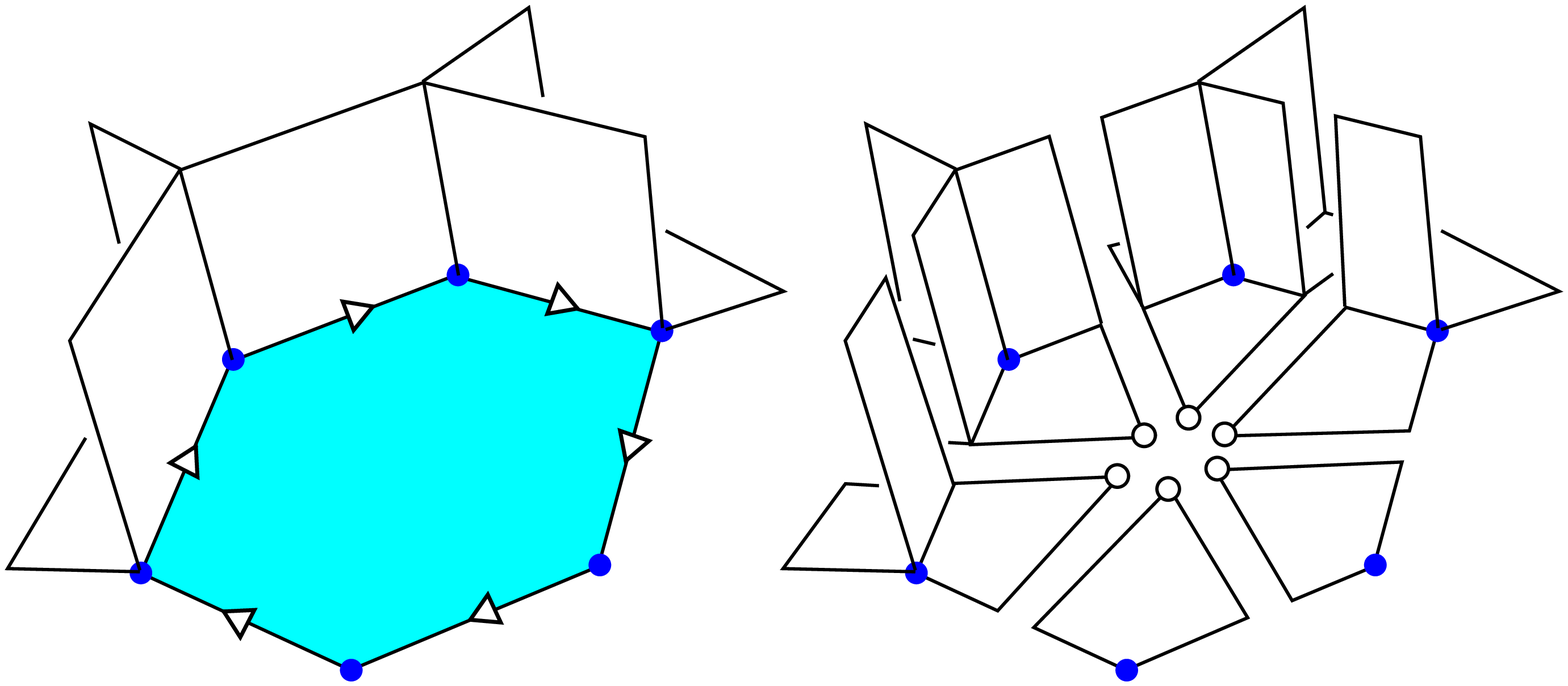}
\nota{A face $f$ of a special complex $X$, with vertices and edges of various types: here $f$ has two vertices of each type $(0,0)$, $(0,1)$, $(0,2)$, and three edges of each type $(1,1)$ and $(1,2)$. We label the oriented edges with some matrices $L_i$ (left) and we break $f$ into star neighborhoods of it vertices (right).}\label{f:face}
\end{figure}

A typical face $f$ of $X$ is shown in Figure \ref{f:face}: it may have vertices and edges of various kinds, and its closure need not to be embedded (it may also be adjacent multiple times to the same edge or vertex). We want to assign labels $L_i$ to the oriented edges (that is, sides) of $f$ as shown in the figure.

\begin{figure}
\centering
\includegraphics[scale=1.4]{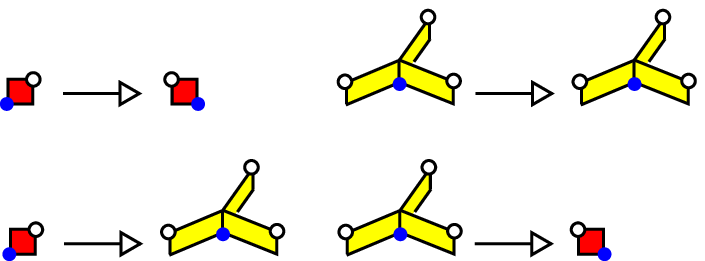}
\nota{Four possible gluings along an oriented edge of type $(1,1)$.}\label{f:blocks_gluing}
\end{figure}

An edge $e_i$ of $f$ can be either of type $(1,1)$ or of type $(1,2)$, and we call it respectively an \emph{interior edge} or a \emph{boundary edge}. A boundary edge is contained in $\partial X$ and connects two vertices $v$ and $v'$ that may be of type $(0,0)$ or $(0,1)$. There are four possible cases, shown in Figure \ref{f:blocks_gluing}. In any case, the fibrations $M_v\to S_v$ and $M_{v'} \to S_{v'}$ get identified along some diffeomorphism $\psi\colon D^2\times S^1 \to D^2\times S^1$ identifying two solid torus facets. As in Section \ref{polygons:section}, we encode this diffeomorphism unambiguously (up to isotopy) via a label 
\[
L_i = \begin{pmatrix} \varepsilon & k \\ 0 & {\varepsilon'} \end{pmatrix}
\]
with $\varepsilon, \varepsilon' = \pm 1$ and $k\in \matZ$. This label is assigned to the side $e_i$ of $f$.

If $e_i$ is an interior edge, it connects two vertices $v$ and $v'$ that may be of type $(0,1)$ or $(0,2)$. The two fibrations $M_v \to S_v$ and $M_{v'} \to S_{v'}$ are now glued along a diffeomorphism $\psi\colon S^1 \times P \to S^1 \times P$ between two facets.

The restriction of $\psi$ to the boundary torus lying above $f$ is a diffeomorphism $S^1 \times S^1 \to S^1 \times S^1$, whose isotopy class is encoded by a matrix $L_i\in {\rm GL}(2,\matZ)$. This is the label that we assign to $e_i$. 

Since the fiber generates the center of $\pi_1(S^1\times P)$, the diffeomorphism $\psi\colon S^1 \times P \to S^1 \times P$ must preserve the fiber (up to reversing the orientation).
Therefore the label $L_i$ has the same nice form as in the previous case: 
\[
L_i = \begin{pmatrix} \varepsilon & k \\ 0 & {\varepsilon'} \end{pmatrix}.
\]
Summing up, a \emph{labeling} of $X$ is simply the assignment of a matrix $\matr \varepsilon k 0 {\varepsilon'}$ to every oriented side $e$ of every face $f$ in $X$. 

We implicitly agree that the orientation reversal of the side $e$ changes the label from $L$ to $L^{-1}$. Note that an interior edge $e$ inherits three labels, one for each incident face, while a boundary edge has only one label.

\subsection{The fibration of $\matCP^2$ over $\Pi_2$}
As an example, we now analyze in detail the labelling on $\Pi_2$ induced by the tropical fibration $\pi:\matCP^2\to \Pi_2$; the answer is depicted in Figure~\ref{f:dots}-(centre), where every unlabeled edge is tacitly assumed to have label $\matr 100{-1}$. This analysis is not necessary for the rest of the paper, so the reader may skip it and jump to Section \ref{adm:subsection}.

Recall that the preimage of all points of type $(k,l)$ with $l\leqslant 1$ is the union of four lines in $\matCP^2$, intersecting in the six points of type $(0,0)$.
Call these lines $l_1,\dots,l_4$, and call $q_1,\dots,q_4$ the four points of type $(0,1)$ corresponding to $l_1,\dots,l_4$ respectively.

Fix an ordered pair $(i,j)$. At the intersection point $p_{ij} = l_i\cap l_j$ we have an identification of a neighborhood $N_{ij}$ of $p_{ij}$ with $D^2\times D^2$ such that $D^2\times \{0\}$ is the intersection of $l_i$ with $N_{ij}$, and $\{0\}\times D^2$ is the intersection of $l_j$ with $N_{ij}$. (The identification is sensitive to swapping $i$ and $j$.) 
Set $q_{ij} = \pi (p_{ij})$.

We fix as above an auxiliary line $r$ going through the points $l_1\cap l_2$ and $l_3\cap l_4$. The line $r$ induces a section of the normal bundles of the four lines, and we use it to fix an identification of all the other facets involved with $S^1 \times P$. With this identification, every internal edge gets a label $\matr {-1}001$, one only needs to check signs by looking at orientations. Using a move that will be described in Proposition \ref{mI:prop}, we can change all these labels with $\matr 100{-1}$.

We need to determine the labels on the external edges. Consider the point $q_{13}$. We are interested in the isotopy class of the section above $l_1$ in the boundary of $N_{13}$: since all lines going through $p_{12}$ are isotopic, the section induced by $r$ on $N_{13}$ is parallel to the curve $S^1\times\{1\}$ in $\partial N_{13}$. Therefore the label on the edge connecting $q_1$ to $q_{13}$ is diagonal, and by looking at the orientations we get $\matr 100{-1}$. Likewise, all edges incident to $q_{14}$, $q_{23}$, and $q_{24}$ are labeled with $\matr 100{-1}$. 

At the point $p_{12}$, the section determined by $r$ on $l_1$ is no longer parallel to $S^1\times\{1\}$ in $\partial N_{12}$.
However, one checks that the section is parallel to the diagonal curve $S^1$ in the corner torus $S^1\times S^1$ in $N_{12}$, and we get $\matr 1{-1}0{-1}$.

Notice that in no case do we need to specify an orientation of the edges, since $\matr100{-1}^2 = \matr 1{-1}0{-1}^2 = I$.

\subsection{Admissibility} \label{adm:subsection}
As in the polygonal case, every fibration is encoded by some (non unique) labeling of $X$, but not every labeling defines a fibration: some simple conditions need to be verified. 

Let $f$ be a face of $X$, with oriented sides $e_1,\ldots, e_n$. Let $v_i$ be the vertex of $f$ adjacent to $e_i$ and $e_{i+1}$. We assign a matrix $J_i$ to $v_i$ as follows:
\begin{itemize}
\item if $v_i$ is of type $(0,0)$, then $J_i = \matr 0110$;
\item if $v_i$ is of type $(0,1)$, then $J_i = \matr 1001$;
\item if $v_i$ is of type $(0,2)$ and is not dotted, then $J_i = \matr 0{1}{1}0$;
\item if $v_i$ is of type $(0,2)$ and is dotted, then $J_i = \matr {-1}011$.
\end{itemize}
Recall that we have fixed two dots at every vertex of type $(0,2)$ as in Figure \ref{f:dots}. Note that in all cases we get $J_i^2 = I$. 

\begin{prop} \label{admissible:prop}
A labeling defines a fibration $M\to X$ if and only if the following hold:
\begin{enumerate}
\item at every oriented interior edge, the three labels of the incident faces are
\[
\begin{pmatrix} \varepsilon & k_1 \\ 0 & {\varepsilon'} \end{pmatrix},
\qquad
\begin{pmatrix} \varepsilon & k_2 \\ 0 & {\varepsilon'} \end{pmatrix},
\qquad
\begin{pmatrix} \varepsilon & k_3 \\ 0 & {\varepsilon'} \end{pmatrix},
\]
for some constants $\varepsilon, \varepsilon' = \pm 1$, with the condition $k_1 + k_2 + k_3 = 0$;
\item at every face $f$ we have
\[
J_nL_n \cdots J_1L_1 = I.
\]
\end{enumerate}
If $\det L_i=-1$ for all $i$, the manifold $M$ is oriented.
\end{prop}
\begin{proof}
At every interior edge we need to build a diffeomorphism $\psi \colon S^1 \times P \to S^1 \times P$, and it is a standard fact in three-dimensional topology that such a diffeomorphism exists if and only if it acts on the boundary tori $S^1 \times S^1$ as specified by condition (1).

Condition (2) is that the monodromy around the central torus must be the identity. The role of $J_i$ is to translate between the two basis of the same corner torus, used by the two adjacent facets. A careful case by case analysis is needed here: 
\begin{itemize}
\item if $v_i$ is of type $(0,0)$, the facets are $S^1 \times D^2$ and $D^2\times S^1$, so $J_i = \matr 0110$; 
\item if $v_i$ is of type $(0,1)$, the facets are $D^2\times S^1$ and $S^1 \times P$, so $J_i = \matr 1001$; 
\item if $v_i$ is of type $(0,2)$, both facets are $S^1 \times P$ and there are two cases:
\begin{itemize}
\item if $v_i$ is not dotted, the factors in $S^1 \times P$ are interchanged as in the case $(0,0)$, so we get $J_i = \matr 0{1}{1}0$, 
\item if $v_i$ is dotted, the boundaries $S^1 \times \partial P$ of the two sections coincide, and we get $J_i = \matr {-1}0{1}1$.
\end{itemize}
\end{itemize}
In the latter case, we have three complex lines $l_1,l_2,r$ passing through a point $p$ and determining three oriented curves $\gamma_1, \gamma_2, \mu$ in the corner torus $S^1\times S^1$. The basis to be compared are $(\gamma_1, \mu)$ and $(\gamma_2, \mu)$ and we have $\mu = \gamma_1 + \gamma_2$, hence $\gamma_2 = \mu-\gamma_1$ and we get $J_i = \matr{-1}011$.
\end{proof}

A labeling on $X$ satisfying the requirements of Proposition \ref{admissible:prop} is \emph{admissible}. If $\det L_i=-1$ then it is \emph{oriented}. An oriented label is either $L = \matr 1k0{-1}$ or $\matr {-1}k01$, and we have called them respectively positive and negative. Note that $L = L^{-1}$ and hence we do not need to orient the edge when assigning it an oriented label. Also in this setting, positive labels are preferable (at least on boundary edges).

\begin{prop}
If all labels are oriented and positive, the nodal surface $S$ is naturally oriented. Every nodal point has positive intersection $+1$.
\end{prop}
\begin{proof}
Same proof as Proposition \ref{positive:prop}, with $P \times D^2$ replacing $D^2\times D^2$.
\end{proof}

We now turn to self-intersection. The nodal surface $S$ is the union of some closed surfaces $S_1 \cup\ldots \cup S_k$ intersecting transversely, such that the abstract resolution of each $S_i$ is connected.

\begin{prop}
If the labels are oriented and positive, and $S_i$ is embedded, then 
\[
S_i\cdot S_i = -\sum_jk_j
\]
as $L_j = \matr 1{k_j}0{-1}$ varies among all labels on edges onto which $S_i$ projects.
\end{prop}
\begin{proof}
Same proof as Proposition \ref{pm:prop}.
\end{proof}

\begin{example}
Consider the two labelings in Figure \ref{f:dots}-(centre) and (right), where every unlabeled edge is tacitly assumed to have label $\matr 100{-1}$. Both labelings are oriented and admissible: the three labels at every interior edges are equal to $\matr 100{-1}$ and hence condition (1) is fulfilled; in the central figure, there are two kinds of faces: the non dotted ones give
\[
J_4L_4\cdots J_1L_1 =  \matr 100{1} \matr 100{-1} \matr 0{1}{1}0 \matr 100{-1}  \matr 100{1} \matr 100{-1}  \matr 0110 \matr 100{-1} = I,
\]
and on the dotted ones we get
\[
J_4L_4\cdots J_1L_1 =  \matr 100{1} \matr 100{-1} \matr {-1}0{1}1 \matr 100{-1}  \matr 100{1} \matr 1{-1}0{-1}  \matr 0110 \matr 1{-1}0{-1} = I.
\]
As seen above, this labeling represents the tropical fibration $\matCP^2 \to \Pi_2$.

On the right figure, we note that there are only two vertices $v$, both of type $(0,1)$, and at every face we have
\[
J_2L_2J_1L_1 = \matr 100{1} \matr 100{-1} \matr 100{1} \matr 100{-1}  = I.
\]
The manifold $M$ here is the double of the basic piece $M_v$ along its boundary. The fiber above $v$ is a $\theta$-shaped graph $\theta$ and $M_v$ is a regular neighborhood of $\theta$, that is a handlebody with one 0--handle and two 1--handles. The double of such a manifold is $M=(S^3\times S^1)\#(S^3\times S^1)$.
\end{example}

\subsection{Moves}
Let $X$ be a special complex equipped with an admissible labeling, defining a fibration $M\to X$. The moves described in Section \ref{moves:subsection} apply also here, and there are more moves that involve vertices of type $(0,1)$ and $(0,2)$ that modify a labeling without affecting the fibration $M\to X$. 

\begin{prop} \label{mI:prop}
Let $v$ be a vertex of $X$, of any type $(0,0)$, $(0,1)$, or $(0,2)$. If we change the signs simultaneously of the labels on all the (two, three, or four) edges incident to $v$, we get a new admissible labeling that encodes the same fibration $M\to X$.
\end{prop}
\begin{proof}
The manifolds $D^2\times D^2$, $D^2 \times P$, and the four-dimensional pair-of-pants $B$ have orientation-preserving self-diffeomorphisms that act like $-I$ on the homologies of all the corner tori in the boundary. 

To see this for $B$, consider $B$ as the complement of some lines in $\matCP^2$ defined by equations with real coefficients. The map $[z_0,z_1,z_2] \mapsto [\bar z_0, \bar z_1, \bar z_2]$ preserves $B$ and acts as required.
\end{proof}

We note in particular that Proposition \ref{blow:prop} is still valid in this context.

\begin{prop}\label{p:special_blowup}
The three moves in Figure \ref{f:moves} transform $M$ into 
\[
M\# \matCP^2, \qquad M\#\CPbar, \qquad M\#(S^2\times S^2)
\]
respectively.
\end{prop}

\begin{rem}
In this section we have dealt only with special complexes, as this simplifies the labelings, but an extension of Propositions \ref{admissible:prop} and \ref{p:special_blowup} to all simple complexes can be done quite easily. In the first proposition, condition (1) is local, and is required also when dealing with nonspecial complexes. Condition (2), on the other hand, is only needed to ensure that the torus fibration on the boundary extends to the interior of the cell; if a connected component of the $(2,2)$--stratum is not a disc, we need to require that the fibration on its boundary extends to the interior. Notice that this extension is not unique in general, hence a labeling in the above sense does not determine a fibration $M\to X$: in order to get uniqueness, we need to specify the monodromy on the boundary as well as its extension. We do not explore this further here.
\end{rem}

\begin{oss}
A 3--manifold decomposing into pieces diffeomorphic to $D^2\times S^1$ and $P\times S^1$ was called a \emph{graph manifold} by Waldhausen \cite{W}: such 3--manifolds are classified and well-understood.
\end{oss}

\section{Fundamental groups} \label{fundamental:groups:section}
In the previous section we have made some effort in defining some labelings that encode all pants decompositions $M\to X$ over a given special complex $X$.
We now use them to prove the following, which is the main result of this paper.

\begin{teo} \label{main:teo}
For every finitely presented group $G$ there is a pants decomposition $M\to X$ with $\pi_1(M) = G$.
\end{teo}

\subsection{Even complexes}
We say that a special complex is \emph{even} if every face is incident to an even number of vertices (counted with multiplicity). Recall that there are three types $(0,0)$, $(0,1)$, and $(0,2)$ of vertices, and each of these must be counted. For instance, the complex $\Pi_2$ is even: every 2--cell is incident to four vertices.

Even complexes are particularly useful here because of the following.
\begin{prop} \label{even:pants:prop}
If $X$ is even, there is a pants decomposition $M\to X$.
\end{prop}
\begin{proof}
We note that every face $f$ in any simple complex contains an even number of vertices of type $(0,1)$. So the evenness hypothesis on $X$ says that the number of vertices of type $(0,0)$ or $(0,2)$ is even for every $f$.

Every vertex $v$ of type $(0,2)$ in $X$ is adjacent to six faces, and we assign dots to two opposite ones arbitrarily. 

We first try to assign trivial labels $L= \matr 1001$ everywhere. Condition (1) of Proposition \ref{admissible:prop} is trivially satisfied, and at every face $f$ we get a product monodromy $J_{2n}L_{2n} \cdots J_1L_1 = J_{2n} \cdots J_1$ that we now compute. 

If there were no dots in $f$, we would get $J_{2n}\cdots J_1 = J^{2k} = I$ with $J=\matr 0110$ and $2k\leqslant 2n$ is the number of vertices of type $(0,0)$ or $(0,2)$. In that case condition (2) would also be satisfied.

If there are some dots, we adjust the labeling so that the above construction still works.
For every maximal string of dotted corners of odd length in a polygonal face, we put a label $\matr 1101$ at the two oriented edges incoming and outcoming the string, both oriented towards the string; e.g.~if there is an isolated dotted corner $v$, the two labels on the edges incoming into $v$ will have label $\matr 1101$, while if there are two connected dotted corners $w,w'$ isolated from all other dotted corners, the label on all the edges incident to $w$ or $w'$ will simply be $\matr 1001$.

This works for the compatibility condition (2), since two consecutive dotted corners contribute with $\matr {-1}011^2 = I$; in an even chain, the product is trivial, while in an odd chain of $2k+1$ dotted vertices we obtain
\[
\matr 1101 \matr {-1}0{1}1^{2k+1} \matr 1{-1}01 = \matr 1101 \matr {-1}0{1}1 \matr 1{-1}01 = \matr  0{1}{1}0.
\]
In either case, after the simplification we are left with a power of $\matr 0110$ for each chain of odd length, and the global monodromy will be trivial for parity reasons (because $X$ is even).

The addition of these labels however may have destroyed condition (1).
Consider an oriented interior edge $e$, that connects either two vertices of type $(0,2)$ or one of type $(0,1)$ and one of type $(0,2)$.

If $e$ is incident to two vertices of type $(0,2)$ there are two possibilities: either the dots are on the same face of $X$ incident to $e$, or they are on different faces.

In the former case, the three labels of $e$ are left unchanged $L=\matr 1001$, and condition (1) is trivially satisfied.
In the latter, two of its three labels have been modified to $\matr 1101$ and $\matr 1{-1}01$ and then condition (1) still holds (note the role of the edge directions).

If $e$ is incident to a vertex $v_2$ of type $(0,2)$ and one $v_1$ of type $(0,1)$, it is incident to three faces, exactly one of which has a dotted corner at $v_2$; denote this face with $f$.
If the label of $e$ as part of $\partial f$ is the trivial label $\matr 1001$, condition (1) is again automatically satisfied.

Suppose now the label of $e$ has been changed. Then condition (1) is violated along $e$, since exactly one label has been modified to $\matr 1101$.
We need to modify the labeling further, and we do so by modifying both the labels at such edges $e$ and on some boundary edges that share a vertex with them.

Consider the set $E_1$ of all external edges $e_1$ with the following property: $e_1$ shares exactly one vertex with an interior edge $e$ such that the label on $e$ on the face $f_1$ that they span is nontrivial (i.e. it is $\matr 1101$).
Let $E_2$ be the set of all external edges $e_2$ with the following property: $e_2$ shares both endpoints with two interior edges $e'$, $e''$, and the labels on $e'$ and $e''$ on the face $f_2$ that they span are both nontrivial (i.e.~they are $\matr 1101$).
By construction, $E_1$ and $E_2$ are disjoint, and so are the associated sets of interior edges.
Also, notice that the faces and edges denoted by $f_1, e$ (respectively, $f_2, e', e''$) are all determined by $e_1$ (resp. $e_2$).

For every edge $e_1$ in $E_1$, we orient it towards $e$, we replace the label of $e_1$ with $\matr 1101$ and let $e$, seen as part of the boundary of $f$, have the trivial label $\matr 1001$. 
For every edge $e_2$ in $E_2$, we replace the two labels on the two associated edges $e'$ and $e''$ (as part of the boundary of $f_2$) by the trivial label $\matr 1001$, and leave the label of $e_2$ unchanged (i.e.~trivial).

It is readily checked that now both conditions (1) and (2) are satisfied.
\end{proof}


There are many even complexes:

\begin{figure}
\labellist
\small\hair 2pt
\pinlabel $v$ at 200 210
\pinlabel $S$ at 270 170
\pinlabel $r_i$ at 100 320
\pinlabel $r_j$ at 170 90
\pinlabel $X$ at 300 20
\pinlabel $X'$ at 700 20
\pinlabel $X''$ at 1100 20
\endlabellist
\centering
\includegraphics[width = 12.5 cm]{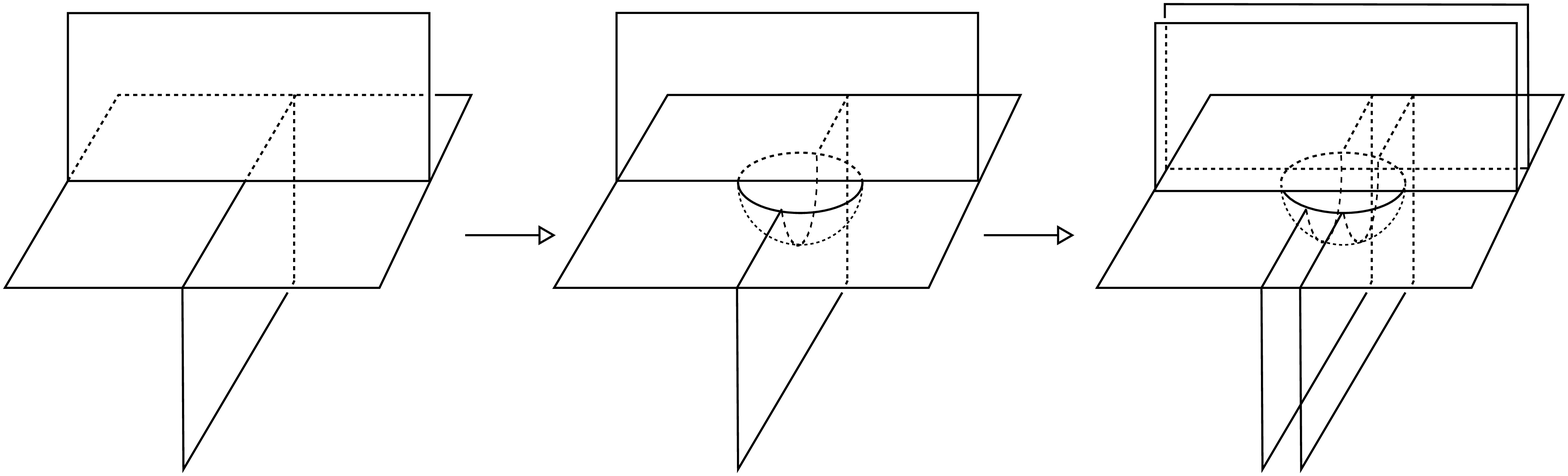}
\nota{How to construct an even complex. The base surface $S$ here is horizontal and the relator faces $r_i$ and $r_j$ are attached vertically.}\label{f:pari}
\end{figure}

\begin{prop} \label{even:group:prop}
Every finitely presented group is the fundamental group of an even complex without boundary.
\end{prop}
\begin{proof}
Every finitely presented group $G =\langle g_1, \ldots, g_k\ |\ r_1, \ldots, r_s \rangle$ is the fundamental group of some special complex $X$ without boundary, constructed by attaching discs to a genus-$k$ surface $S$. To see this, first attach discs to $S$ to transform the fundamental group of $S$ into a free group $F_k$ with $k$ generators (for instance, you may take the meridians of a handlebody with boundary $S$). Then attach discs on $S$ along $s$ generic curves that represent the relators $r_1,\ldots, r_s$ in a sufficiently complicated way, so that $S$ is cut into polygons by them (add a trivial relator $r_1$ in case there are none).

We now modify $X$ to an even complex $X''$ with the same fundamental group $G$. The modification is depicted in Figure \ref{f:pari} and consists of two steps: the first is a local modification at every vertex $v$ of $S$, where two relator faces $r_i$ and $r_j$ intersect. Note that every face in $X'$ is even, except the new small triangles created by the move. Then we double each relator $r_i$ as shown in the figure (that is, for every $i=1,\dots, s$ we attach two parallel discs), Now triangles are transformed into squares: the final polyhedron $X''$ is even and has the same fundamental group $G$ of $X$ and $X'$.
\end{proof}

\subsection{Fundamental group}
How can we calculate the fundamental group of $M$ by looking at the fibration $M\to X$? We answer this question in some cases. We start by showing that in  dimension $4$ any facet of the compact pair of pants carries the fundamental group of the whole block (in contrast with dimension 2).

\begin{lemma} \label{surjective:lemma}
Let $B$ be the compact $4$--dimensional pair of pants and $F\cong P\times S^1$ be any of its four facets. The map $\pi_1(F) \to \pi_1(B)$ induced by inclusion is surjective.
\end{lemma}
\begin{proof}
Recall that $B$ is $\matCP^2$ minus the open regular neighborhood of four lines $l_1,l_2,l_3,l_4$. Let $F$ correspond to $l_4$. Using the Salvetti complex \cite{S} we see that $\pi_1(B) \cong \matZ^3$ is generated by three loops turning around any three of these lines, say $l_1, l_2, l_3$. These loops can be homotoped inside $F \cong P\times S^1$, where they correspond to three meridians on the boundary tori.
\end{proof}

Let $X^1$ denote the 1--stratum of $X$, that is the set of all non-smooth points of $X$.

\begin{prop}
Let $M\to X$ be a pants decomposition. The induced map $\pi_1(M) \to \pi_1(X)$ on fundamental groups is surjective. It is also injective, provided the following holds:
\begin{itemize}
\item $X$ is not a surface, 
\item every connected component of $X^1 \setminus \partial X$ is incident to a vertex in $\partial X$ whose fiber is contained in a (possibly immersed) spherical component of the nodal surface.
\end{itemize}
\end{prop}
\begin{proof}
The map $\pi_1(M) \to \pi_1(X)$ is surjective because all fibers are connected and arcs lift from $X$ to $M$. 

Let $F_x = \pi^{-1}(x)$ be the fiber of $x$ and let $G_x$ be the image of the map $\pi_1(F_x) \to \pi_1(M)$ induced by inclusion (with some basepoint in $F_x$).
It is easy to prove that if $G_x$ is trivial for every $x\in X$, then $\pi_1(M) \to \pi_1(X)$ is an isomorphism. We now prove that the additional assumptions listed above force all groups $G_x$ to be trivial.

We use the term \emph{connected stratum} to denote a connected component of some $(k,l)$--stratum of $X$. If $G_x$ is trivial for some $x$, then $G_{x'}$ is trivial for all points $x'$ lying in the same connected stratum of $x$ and we say that the connected stratum is \emph{trivial}. We now show that the triviality propagates along incident connected strata in most (but not all!) cases. Let $s$ and $t$ we two incident connected strata, that is such that either $s\subset \bar t$ or $t\subset \bar s$. Suppose that $s$ is trivial. We claim that, if any of the following conditions holds, then $t$ is also trivial.
\begin{enumerate}
\item $\dim t > \dim s$; 
\item $t \subset \partial X$, $s\not\subset \partial X$, and $\dim t = \dim s-1$; 
\item $t$ is a vertex of type $(0,2)$ and $s$ is an edge of type $(1,2)$. 
\end{enumerate}
To prove the claim, pick $x\in s$ and $y \in t$; by assumption, $G_x$ is trivial.

\begin{enumerate}
\item We have $s\subset \bar t$ and the fiber $F_y$ can be isotoped to $F_{y'}$ where $y'$ is close to $x$, so $F_{y'}$ lies in a regular neighborhood of $F_x$, therefore $G_y$ is naturally a subgroup of $G_x$, hence trivial.

\item In particular case $F_x \cong F_y\times S^1$ and $F_y$ can be isotoped inside $F_x$.

\item It follows from Lemma~\ref{surjective:lemma}.
\end{enumerate}

By assumption every connected component $C$ of $X^1 \setminus \partial X$ is incident to a vertex $v$ of type $(0,1)$ in $\partial X$, whose fiber $F_v$ is contained in a sphere: therefore $G_v$ is trivial. By property (2) the edge of type (1,2) adjacent to $v$ is also trivial, and we can use (1) and (3) to propagate the triviality along all the connected strata of $C$.

Since $X$ is not a surface, every 2--dimensional connected stratum of $X$ is incident to $X^1\setminus \partial X$, and is hence trivial by property (1). Finally, the triviality extends to the rest of $\partial X$ by (2).
\end{proof}

\begin{cor} \label{isomorfismo:cor}
Let $M\to X$ be a pants decomposition. If $X^1 \setminus \partial X$ is connected, $\partial X \neq \varnothing$, and the nodal surface consists of (possibly immersed) spheres, the map $\pi_1(M) \to \pi_1(X)$ is an isomorphism.
\end{cor}

The homomorphism $\pi_1(M) \to \pi_1(X)$ may not be injective in general: Figure \ref{f:dots}-(right) shows a fibration $M\to X$ with $\pi_1(M) = \matZ * \matZ$ and $\pi_1(X) = \{e\}$.

\subsection{Proof of the main theorem}
We can finally prove the main result of this paper, that is Theorem \ref{main:teo}.

\begin{figure}
\labellist
\small\hair 2pt
\pinlabel $X$ at 100 10
\pinlabel $e$ at 73 75
\pinlabel $X'$ at 370 10
\endlabellist
\centering
\includegraphics[scale=.7]{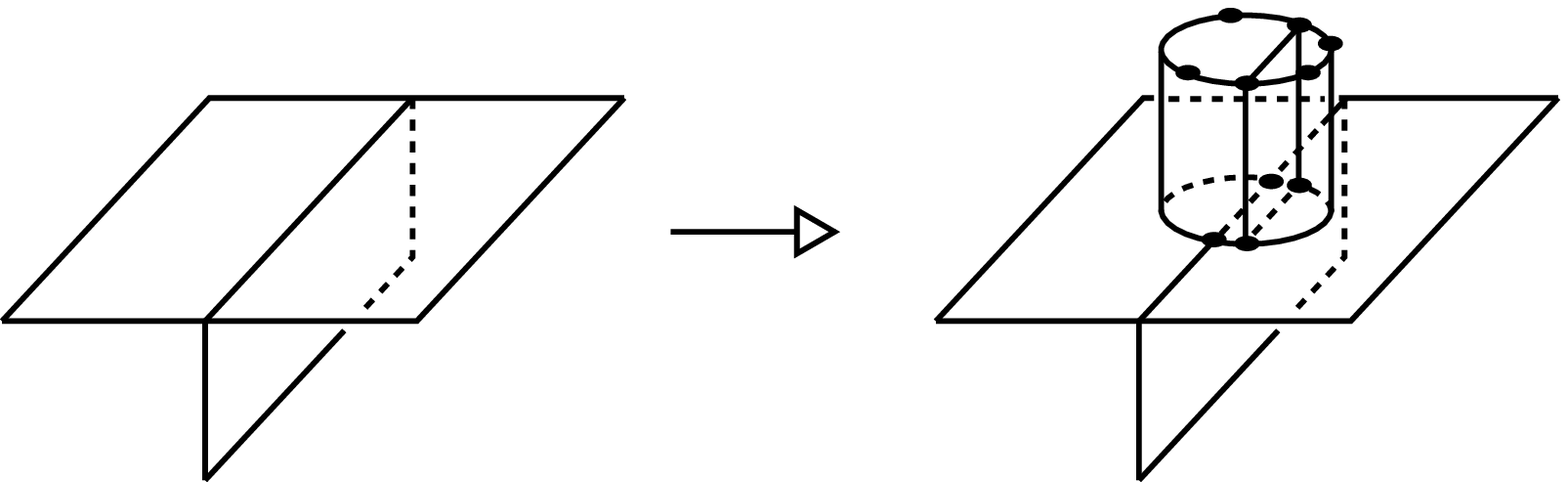}
\nota{How to create some boundary on an even complex, preserving evenness and the fundamental group. The complex $X'$ is constructed by attaching a product $\theta \times [0,1]$ to $X$ along $\theta \times 0$ as shown, where $\theta$ is a $\theta$-shaped graph. The boundary $\partial X' = \theta \times 1$ contains four vertices of type $(0,0)$ and two of type $(0,1)$: all dotted points are vertices of some type.}\label{f:create_boundary}
\end{figure}

\begin{proof}[Proof of Theorem~\ref{main:teo}]
For every finitely presented group $G$ there is an even special complex $X$ without boundary and with $\pi_1(X) = G$ by Proposition \ref{even:group:prop}. We modify slightly $X$ to a complex $X'$ with non-empty boundary, by choosing an arbitrary edge $e$ and modifying $X$ near $e$ as shown in Figure~\ref{f:create_boundary}. 

We have $\pi_1(X) = \pi_1(X')$ and $X'$ is still even. Note that $\partial X'$ is a $\theta$-shaped graph with two vertices of type $(0,1)$ and also four vertices of type $(0,0)$, indicated in the picture. Note also that $X^1$ is connected because $X$ is special without boundary, and hence $(X')^1 \setminus \partial X'$ is also connected.

By Proposition \ref{even:pants:prop} there is a pants decomposition $M \to X'$. By looking at $\partial X'$ we see that the nodal curve consists of three spheres. Corollary  \ref{isomorfismo:cor} hence applies and gives $\pi_1(M) = \pi_1(X') = G$.
\end{proof}


\begin{thebibliography}{1}

\bibitem{CoThu}
Fran\c{c}ois Costantino and Dylan~P. Thurston, \emph{3--manifolds efficiently
  bound 4--manifolds}, J. Topol. \textbf{1} (2008), no.~3, 703--745.

\bibitem{FY}
Stephan Fischli and David Yavin, \emph{Which 4--manifolds are toric
  varieties?}, Math. Z. \textbf{215} (1994), no.~2, 179--185.

\bibitem{Mat} Sergej Matveev, ``Algorithmic topology and classification of 3-manifolds'', second edition, Algorithms and Computation in Math. \textbf{9}, Springer, Berlin (2007).

\bibitem{M}
Grigory Mikhalkin, \emph{Decomposition into pairs-of-pants for complex
  algebraic hypersurfaces}, Topology \textbf{43} (2004), no.~5, 1035--1065.

\bibitem{PR}
Mikael Passare and Hans Rullg{\aa}rd, \emph{Amoebas, monge-amp\`ere measures
  and triangulations of the newton polytope}, Duke Math. J. \textbf{121}
  (2004), no.~3, 481--507.

\bibitem{S}
Mario Salvetti, \emph{Topology of the complement of real hyperplanes in
  $\matc^n$}, Invent. Math. \textbf{88} (1987), no.~3, 603--618.
  
\bibitem{Tu} 
Vladimir Turaev, \emph{Quantum invariants of knots and 3-manifolds}, De Gruyter Studies in Mathematics, \textbf{18}, 2010.

\bibitem{W}
Friedhelm Waldhausen, \emph{Eine Klasse von 3--dimensionalen
  Mannigfaltigkeiten. I}, Invent. Math. \textbf{3} (1967), 308--333.

\end{thebibliography}
\end{document}